\documentclass[aip,
amsmath,amssymb,
preprint,%
numerical,
]{revtex4-1}
\usepackage{graphicx}% Include figure files
\usepackage{dcolumn}% Align table columns on decimal point
\usepackage{bm}% bold math
\usepackage{hyperref}
\usepackage[utf8]{inputenc}
\usepackage[T1]{fontenc}
\usepackage{mathptmx}
\usepackage{etoolbox}
\newcommand{\RNum}[1]{\uppercase\expandafter{\romannumeral #1\relax}}
\usepackage{array,xcolor,colortbl}
\usepackage{multirow}
\makeatletter
\def\@email#1#2{%
	\endgroup
	\patchcmd{\titleblock@produce}
	{\frontmatter@RRAPformat}
	{\frontmatter@RRAPformat{\produce@RRAP{*#1\href{mailto:#2}{#2}}}\frontmatter@RRAPformat}
	{}{}
}%
\makeatother
\begin{document}

	\title{Effect of rotation at the bio-convective instability in a suspension of isotropic scattering phototactic algae.}
	% Force line breaks with \\

		% Force line breaks with \\
	\author{S. K. Rajput}
	\altaffiliation[Email: ]{shubh.iitj@gmail.com}
	%\author{M. K. Panda}%
	%\email{mkpanda@iiitdmj.ac.in}
	%\email{shubh.iiitj@gmail.com}
	\affiliation{ 
		Department of Mathematics, PDPM Indian Institute of Information Technology Design and Manufacturing,
		Jabalpur 482005, India.%\\This line break forced with \textbackslash\textbackslash
	}%
	%\homepage{http://www.Second.institution.edu/~Charlie.Author.}
	%\affiliation{%
		%	Second institution and/or address%\\This line break forced% with \\
		%}%
	

	\begin{abstract}
		
	 The impact of rotation on the bio-convective instability in a phototactic algal suspension is investigated through a linear stability analysis in this article. The study primarily focuses on the suspension's phototactic behavior, which refers to the movement of microorganisms towards or away from light. In conditions of low light, microorganisms exhibit positive phototaxis, moving towards the light source, whereas in intense light, they display negative phototaxis, moving away from the light source. The study specifically examines a suspension that is illuminated by collimated flux, with a constant radiation intensity applied to the top surface. The stability analysis is conducted through the implementation of linear perturbation theory, which enables the prediction of both stationary and oscillatory characteristics of the bio-convective instability at the onset of biocnvection. It is observed that rotation plays a significant stabilizing role in the system.

%		
%		\textbf{ In the natural environment, the sunlight gets scattered by some atmospheric effects (clouds). Therefore, the light reaching the earth's surface divides into two parts collimated and diffuse. So, there is a need to investigate the effect of both types of light (flux). In the Industries, sometimes need a controlled convection which can be achieved to consider rigid boundaries. In this article, we investigate the suspension stability whose upper and top boundaries are rigid no-slip, and which is illuminated by both types of light. Due to the rigidness of the walls, there is a need for a higher Rayleigh number to the onset of phototactic bioconvection as compared to the case of a stress-free top surface. In presence of the diffuse irradiation, the suspension becomes more stable and the oscillatory behavior of the solution is frequently observed.  }
		
	\end{abstract}

	\maketitle

	\section{INTRODUCTION}
	
	Bioconvection, a fascinating phenomenon, encompasses the convective motion observed in fluid containing self-propelled motile microorganisms like algae and bacteria at a macroscopic level ~\cite{20platt1961,21pedley1992,22hill2005,23bees2020,24javadi2020}. These microorganisms exhibit a tendency to move upwards on average due to their higher density compared to the surrounding medium, usually water. When the microorganisms cease their movement, the distinct pattern formation in bioconvection disappears. However, it's important to note that pattern formation in bioconvection is not solely dependent on upswimming or higher microorganism density. Instead, microorganisms respond to various environmental stimuli, referred to as "taxes," by altering their swimming direction. Examples of such taxes include gravitaxis (response to gravitational acceleration), chemotaxis (response to chemicals), phototaxis (response to illumination intensity), and gyrotaxis (generated by two torques due to gravity and local shear flow). In this article, we specifically focus on examining the impacts of phototaxis.
	
	Experimental studies have provided insights into the influence of  collimated flux on bioconvective patterns~\cite{1wager1911,2kitsunezaki2007}. The intensity of light plays a crucial role in shaping the stable patterns observed in suspensions of motile microorganisms in well-stirred cultures. Bright light can either disrupt the existing patterns or prevent pattern formation altogether. The size, shape, structure, and symmetry of the patterns can be significantly affected by the light intensity~\cite{3kessler1985,4williams2011,5kessler1989}.
	The deformation of bioconvection patterns can be explained by considering the interaction between light intensity and the phototactic behavior of the cells. When the illumination intensity $G$ is below or above a critical value $G_c$, the cells exhibit positive or negative phototaxis, respectively, to optimize photosynthesis or avoid photodamage. Consequently, the cells tend to accumulate in regions with optimal light conditions where $G = G_c$. Another factor influencing pattern modifications is the absorption and scattering of light by the cells. Absorption results in a decrease in light intensity along the incident path, while scattering causes light to deviate from the incident path, leading to intensity fluctuations. These fluctuations can cause the light to scatter to different points, resulting in localized increases or decreases in intensity~\cite{7ghorai2010}.

	The study utilizes the phototaxis model proposed by S. Kumar (referred to as SK), the Navier-Stokes equations and a microorganism conservation equation are employed. The suspension of phototactic microorganisms is assumed to be rotating around the z-axis at a constant angular velocity and is illuminated from above with collimated irradiation. As many motile algae rely on photosynthesis for food production, they exhibit pronounced phototaxis. To accurately capture their behavior, it is crucial to analyze the phototaxis model within a rotating medium.
	However, the original SK model neglected the scattering effect caused by the microorganisms. As a result, the microorganisms only received light directly from the source located directly above them, leading to a purely vertical swimming direction. The swimming velocity was determined based on the light intensity reaching each cell. In our modified model, we incorporate the scattering effect caused by the cells. This means that the microorganisms now interact with light through scattering, resulting in a more complex swimming behavior. As a consequence, the swimming direction of the microorganisms includes a horizontal component in addition to the vertical component. By considering scattering, our model allows for a more accurate representation of the swimming behavior of microorganisms in response to light.
	
	In a suspension of finite depth that is illuminated from above, the basic steady state is characterized by a balance between phototaxis (movement of microorganisms in response to light) due to absorption and scattering, and diffusion caused by the random swimming motions of the cells. This balance leads to the formation of a concentrated layer of microorganisms, known as the sublayer, which is horizontally oriented. The position of the sublayer within the suspension depends on the critical light intensity ($G_c$).
	
	If the light intensity across the suspension is lower (higher) than $G_c$, the sublayer forms at the top (bottom) of the chamber. If $G_c$ lies between the maximum and minimum intensities across the suspension, a sublayer forms within the interior of the chamber. In relation to the sublayer, the region below (above) it is gravitationally unstable (stable). Consequently, if the fluid layer becomes unstable, the fluid motions in the unstable layer penetrate into the stable layer. This phenomenon is known as penetrative convection and is observed in various convection problems~\cite{9straughan1993,10ghorai2005,11panda2016}.

Phototactic bio-convection (PBC) has been the subject of extensive research in the scientific literature. Vincent and Hill~\cite{12vincent1996} simulated proprietary PBC and identified the negative buoyancy of cells as the primary factor affecting the suspension. They employed the Beer-Lambert law to model light intensity and conducted linear stability analysis. Ghorai and Hill~\cite{10ghorai2005} quantified PBC in two dimensions and addressed the deficiencies in Vincent and Hill's equilibrium solution. Ghorai et al.~\cite{7ghorai2010} studied the impact of light scattering and discovered bimodal equilibrium configurations in nearly purely scattering suspensions. Ghorai and Panda~\cite{14panda2013} investigated the initiation of bioconvection in anisotropic scattering solutions of phototactic algae, while Panda and Singh~\cite{11panda2016} numerically explored the linear stability of PBC using Vincent and Hill's continuum model in a two-dimensional setting. Previous studies did not account for the effects of diffuse irradiation. Panda et al.~\cite{15panda2016} proposed a model that examined how diffuse radiation influenced isotropic scattering algal suspensions, revealing that diffuse irradiation significantly stabilizes such suspensions. Panda~\cite{8panda2020} further investigated the impact of forward anisotropic scattering on the onset of PBC, considering both diffuse and collimated irradiation. Panda $et$ $al$.~\cite{16panda2022}, in our study, analyzed the impact of oblique collimated irradiation on algal suspensions and checked the linear stability of the suspension for different angles of incidence. Kumar~\cite{17kumar2022} examined the effect of oblique collimated irradiation on isotropic scattering algal suspensions, revealing the existence of both stationary and overstable solutions within certain parameter ranges. Kumar~\cite{39kumar2023} explored the effect of collimated irradiation on algae suspensions, assuming rigid vertical walls and observing a stabilizing effect on the suspension due to these walls. Recently, Kumar~\cite{40kumar2023} investigated the effect of rotation on the phototactic bioconvection in a non scattering medium and foound a significant stabilizing imapact due to rotation. More recently, Panda and Rajput~\cite{41rajput2023} investigate the combined impact of oblique and diffuse irradiation on the onset of bioconvection.
 Building on Kumar's recent work, we investigates the influence of a rotation on the onset of instability in an isotropic scattering algal suspension that is illuminated by vertical collimated flux.
	
The article is divided into three primary sections. Firstly, the authors propose a mathematical model to study phototaxis and calculate the equilibrium state. In the next section, they analyze the linear stability of the equilibrium state using perturbation theory and solve it numerically using the Newton-Raphson-Kantorovich method, obtaining neutral curves. Finally, the authors discuss and interpret the results based on these neutral curves.
	
	\section{MATHEMATICAL FORMULATION}
	
	In this study, the focus is on the movement of a dilute phototactic algal suspension within a confined space between two parallel horizontal boundaries. The suspension is subjected to illumination from above, with collimated flux. The depth of the suspension is fixed at $H$, and the boundaries are impermeable. The light intensity at any specific location $\boldsymbol{x}$ within the suspension, in a particular unit direction $\boldsymbol{r}$, is denoted by $L(\boldsymbol{x}, \boldsymbol{r})$.
	
	\subsection{PHOTOTAXIS IN ALGAE SUSPENSION WITH ABSORPTION AND SCATTERING EFFECTS}
	The radiative transfer equation (RTE) provides a mathematical framework for describing the behavior of radiation in a medium, including an algal suspension with absorption and scattering. It allows us to calculate the light intensity at a particular location $\boldsymbol{x}$ and in a specific direction $\boldsymbol{r}$ within the suspension. The RTE equation, in its general form, can be written as:
	\begin{equation}\label{1}
		\boldsymbol{ r}\cdot\nabla L(\boldsymbol{x},\boldsymbol{r})+(a+\sigma)L(\boldsymbol{x},\boldsymbol{r})=\frac{\sigma}{4\pi}\int_{0}^{4\pi}L(\boldsymbol{x},\boldsymbol{r'})P(\boldsymbol{r},\boldsymbol{r'})d\Omega',
	\end{equation}
The absorption coefficient, denoted as $a$, and the scattering coefficient, represented by $\sigma$, are parameters used to describe the interaction of light with the algal suspension. The scattering phase function, $P(\boldsymbol{r},\boldsymbol{r'})$, characterizes the angular distribution of scattered light. In the specific case of this study (Panda et al.\cite{15panda2016}), the algae cells assume that light scattering occurs in an isotropic manner, simplifying the scattering phase function to a constant value of 1. This simplification is adopted for the sake of simplicity in subsequent calculations, as mentioned in the work of Panda et al.\cite{15panda2016}.

	\begin{figure}[!h]
		\centering
		\includegraphics[width=14cm]{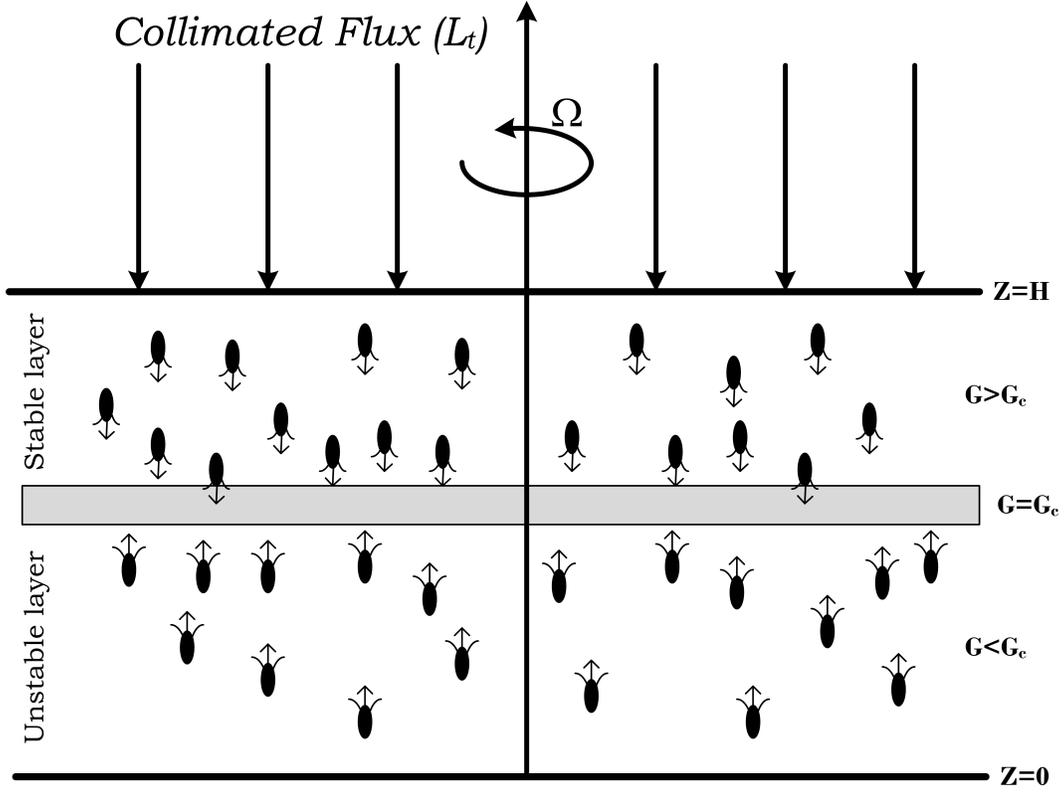}
		\caption{\footnotesize{Axial representation of the problem.}}
		\label{fig1}
	\end{figure}

	Consider the suspension's upper boundary is diffusive and hence, the light intensity at the suspension's top can be represented as
	\begin{equation*}
		L(\boldsymbol{x}_b,\boldsymbol{r})=L_t\delta(\boldsymbol{r}-\boldsymbol{r_0}), 
	\end{equation*}
	where $\boldsymbol{x}_b=(x,y,H)$ is the location on the top boundary surface. Here, $L_t$ is the magnitude of collimated flux. Consider $a=\alpha n(\boldsymbol{x})$ and $\sigma=\beta n(\boldsymbol{x})$. With these substitutions, Eq.~(\ref{1}) can be rewritten as
	\begin{equation}\label{2}
		\boldsymbol{ r}\cdot\nabla L(\boldsymbol{x},\boldsymbol{r})+(\alpha+\beta)nL(\boldsymbol{x},\boldsymbol{r})=\frac{\beta n}{4\pi}\int_{0}^{4\pi}L(\boldsymbol{x},\boldsymbol{r'})d\Omega'.
	\end{equation}
	The value of the total intensity at a given point $\boldsymbol{x}$ is determined by 
	\begin{equation*}
		G(\boldsymbol{x})=\int_0^{4\pi}L(\boldsymbol{x},\boldsymbol{r})d\Omega,
	\end{equation*}
	and similarly, the radiative heat flux is given by 
	\begin{equation}\label{3}
		\boldsymbol{q}(\boldsymbol{x})=\int_0^{4\pi}L(\boldsymbol{x},\boldsymbol{r})\boldsymbol{r}d\Omega.
	\end{equation}
	We assume that the cells and fluid flow at the same speed. Therefore, we can define the cells mean swimming velocity as
	\begin{equation*}
		\boldsymbol{W}_c=W_c<\boldsymbol{p}>,
	\end{equation*}
	here, $W_c$ denotes the cells mean swimming speed average swimming speed, and $<\boldsymbol{p}>$ denotes the mean direction of the cell's swimming which is determined by using the following equation
	\begin{equation}\label{4}
		<\boldsymbol{p}>=-M(G)\frac{\boldsymbol{q}}{|\boldsymbol{q}|}.
	\end{equation}
	Here, taxis response function (taxis function) $M(G)$ describes how algae cells react to light and take a mathematical form as 
	\begin{equation*}
		M(G)=\left\{\begin{array}{ll}\geq 0, & \mbox{ } G(\boldsymbol{x})\leq G_{c}, \\
			< 0, & \mbox{ }G(\boldsymbol{x})>G_{c}.  \end{array}\right. 
	\end{equation*}
	When the light intensity reaches a critical value ($G = G_c$), the microorganisms exhibit zero mean swimming direction. The form of the taxis function varies depending on the species of the microorganisms.~\cite{12vincent1996} 
	
	\subsection{GOVERNING EQUATIONS}
	
	Assume that a continuous distribution is used to model cell population, as has been done in previous studies. In an incompressible dilute algal suspension, each algal cell has a volume $V$ and a density of $\rho+\Delta\rho$, where $\rho$ is the density of water. In this model, it is assumed that the all physical properties of the fluid are constant except the buoyancy force. The governing equations for the system are given by the following equations\\
	1. Continuity equation
	\begin{equation}\label{5}
		\boldsymbol{\nabla}\cdot \boldsymbol{u}=0,
	\end{equation}
	where $\boldsymbol{u}$ is fluid velocity.\\
	2. In the rotating medium, momentum equation under Boussinsque approximation
	\begin{equation}\label{6}
		\rho\left(\frac{\partial \boldsymbol{u}}{\partial t}+(\boldsymbol{u}\cdot\nabla )\boldsymbol{u}+2\boldsymbol{\Omega}\times \boldsymbol{u}\right)=-\boldsymbol{\nabla} P_e+\mu\nabla^2\boldsymbol{u}-nVg\Delta\rho\hat{\boldsymbol{z}},
	\end{equation}
	where $g$ is the gravitational acceleration due to gravity, $\boldsymbol{\Omega}=\Omega\hat{\boldsymbol{z}}$ is the angular velocity, $P_e$ is the excess pressure above hydrostatic, and $\mu$ is the viscosity of the fluid.\\
	3. Cell conservation equation
	\begin{equation}\label{7}
		\frac{\partial n}{\partial t}=-\boldsymbol{\nabla}\cdot \boldsymbol{F_1}=\nabla\cdot[n\boldsymbol{u}+nW_c<\boldsymbol{p}>-\boldsymbol{D}\boldsymbol{\nabla} n].
	\end{equation}
	Here, two key assumptions are made. First, the microorganisms are purely phototactic, and second, $\boldsymbol{ D}=DI$. With the help of these two assumptions, we can remove the Fokker-Plank equation from the governing system.~\cite{15panda2016}

	In this model, the lower horizontal boundary is assumed to be rigid and upper horizontal boundary is assumed to be stress-free. Therefore, the boundary conditions can be expressed as
	\begin{equation}\label{8}
		\boldsymbol{u}\cdot\hat{\boldsymbol{z}}=0\qquad at\quad z=0,H,
	\end{equation}
	\begin{equation}\label{9}
		\boldsymbol{F_1}\cdot\hat{\boldsymbol{z}}=0\qquad at\quad z=0,H,
	\end{equation}
	\begin{equation}\label{10}
	\boldsymbol{u}\times\hat{\boldsymbol{z}}=0\qquad at\quad z=0,
	\end{equation}
	\begin{equation}\label{11}
	\frac{\partial^2}{\partial z^2}(\boldsymbol{u}\cdot\hat{\boldsymbol{z}})=0\qquad at\quad z=H.
	\end{equation}
	
	The top boundary is assumed to be exposed to collimated direct irradiation, then the boundary conditions for intensities are as follows
	\begin{subequations}
		\begin{equation}\label{12a}
			at~~~ z=H,~~~~~~~	L(x, y, z , \theta, \phi)=L_t\delta(\boldsymbol{s}-\boldsymbol{s_0}),~~~ where~~~ (\pi/2\leq\theta\leq\pi),
		\end{equation}
		\begin{equation}\label{12b}
			at~~~ z=0,~~~~~~~ 	L(x, y, z , \theta, \phi) =0,~~~ where~~~ (0\leq\theta\leq\pi/2).
		\end{equation}
	\end{subequations}

	\subsection{DIMENSIONLESS EQUATIONS}
	The equations that govern the system are transformed into a dimensionless form by selecting suitable scales for length ($H$), time ($H^2/D$), velocity ($D/H$), pressure ($\mu D/H^2$), and concentration ($\tilde{n}$). This is done to simplify the equations and make them easier to solve. The resulting dimensionless equations are presented below
	\begin{equation}\label{13}
		\boldsymbol{\nabla}\cdot\boldsymbol{u}=0,
	\end{equation}
	\begin{equation}\label{14}
		S_c^{-1}\left(\frac{\partial \boldsymbol{u}}{\partial t}+(\boldsymbol{u}\cdot\nabla )\boldsymbol{u}\right)+\sqrt{T_a}(\hat{z}\times\boldsymbol{u})=-\nabla P_{e}-Rn\hat{\boldsymbol{z}}+\nabla^{2}\boldsymbol{u},
	\end{equation}
	\begin{equation}\label{15}
		\frac{\partial{n}}{\partial{t}}=-\boldsymbol{\nabla}\cdot\boldsymbol{F_1}=-{\boldsymbol{\nabla}}\cdot[\boldsymbol{n{\boldsymbol{u}}+nV_{c}<{\boldsymbol{p}}>-{\boldsymbol{\nabla}}n.}]
	\end{equation}
	
	In the above equations, $S_c^{-1}=\nu/D$ represents the Schmidt number, $V_c$ denotes the dimensionless swimming speed as $V_c=W_cH/D$, $R=\tilde{n}V g\Delta{\rho}H^{3}/\nu\rho{D}$ is Rayleigh number, and $T_a=4\Omega^2H^4/\nu^2$ is the Taylor number.
	
	After non-dimensionalization, the boundary conditions are expressed as
	\begin{equation}\label{16}
	\boldsymbol{u}\cdot\hat{\boldsymbol{z}}=0\qquad at\quad z=0,1,
	\end{equation}
	\begin{equation}\label{17}
	\boldsymbol{F_1}\cdot\hat{\boldsymbol{z}}=0\qquad at\quad z=0,1,
	\end{equation}
	\begin{equation}\label{18}
	\boldsymbol{u}\times\hat{\boldsymbol{z}}=0\qquad at\quad z=0,
	\end{equation}
	\begin{equation}\label{19}
	\frac{\partial^2}{\partial z^2}(\boldsymbol{u}\cdot\hat{\boldsymbol{z}})=0\qquad at\quad z=1.
	\end{equation}

	The RTE in dimensionless form is
	\begin{equation}\label{20}
		\boldsymbol{r}\cdot\nabla L(\boldsymbol{x},\boldsymbol{r})+\kappa nL(\boldsymbol{x},\boldsymbol{r})=\frac{n\sigma_s }{4\pi}\int_{0}^{4\pi}L(\boldsymbol{x},\boldsymbol{r'})d\Omega',
	\end{equation}
	where $\kappa=(\alpha+\beta)\tilde{n}H$ is the dimensionless absorption coefficient and $\sigma_s=\beta\tilde{n}H$ is the dimensionless scattering coefficient. The scattering albedo $\omega=\sigma_s/\kappa$ measures the scattering efficiency of microorganisms. The Eq.~(\ref{16}) can be alternatively written in terms of the scattering albedo $\omega$ as follows
	\begin{equation}\label{21}
		\boldsymbol{r}\cdot\nabla L(\boldsymbol{x},\boldsymbol{r})+\kappa nL(\boldsymbol{x},\boldsymbol{r})=\frac{\omega\kappa n}{4\pi}\int_{0}^{4\pi}L(\boldsymbol{x},\boldsymbol{r'})d\Omega',
	\end{equation}
	where $\omega\in[0,1]$. A value of $\omega=1$ indicates a purely scattering medium, while $\omega=0$ corresponds to a purely absorbing medium. The RTE can also be written as follows
	
	\begin{equation}\label{22}
		\xi\frac{dL}{dx}+\eta\frac{dL}{dy}+\nu\frac{dL}{dz}+\kappa nL(\boldsymbol{x},\boldsymbol{r})=\frac{\omega\kappa n}{4\pi}\int_{0}^{4\pi}L(\boldsymbol{x},\boldsymbol{r'})d\Omega',
	\end{equation}
	where $\xi,\eta$ and $\nu$ are the direction cosines in the direction $\boldsymbol{r}$. In dimensionless form, the boundary conditions for the intensity are
	\begin{subequations}
		\begin{equation}\label{23a}
			at~~~ z=1,~~~~~~~	L(x, y, z, \theta, \phi)=L_t\delta(\boldsymbol{r}-\boldsymbol{r_0}),~~~ where~~~ (\pi/2\leq\theta\leq\pi),
		\end{equation}
		\begin{equation}\label{23b}
			at~~~ z=0,~~~~~~~	L(x, y, z, \theta, \phi) =0,~~~ where~~~ (0\leq\theta\leq\pi/2). 
		\end{equation}
	\end{subequations}

	\section{THE STEADY SOLUTION}
	
	The equations $(\ref{13})-(\ref{15})$ and $(\ref{22})$ have an equilibrium solution that can be described by the following equation
	
	\begin{equation}\label{24}
		\boldsymbol{u}=0,~~~\zeta_s=\nabla\times u=0,~~~n=n_s(z)\quad and\quad  L=L_s(z,\theta),
	\end{equation}
	where, $\zeta_s$, is the vorticity vector at the basic steady state.
	Thus, in the equilibrium state, the total intensity $G_s$ and radiative heat flux $\boldsymbol{q}_s$ can be expressed as follows
	\begin{equation*}
		G_s=\int_0^{4\pi}L_s(z,\theta)d\Omega,\quad 
		\boldsymbol{q}_s=\int_0^{4\pi}L_s(z,\theta)\boldsymbol{s}d\Omega,
	\end{equation*}
	and the governing equation for $L_s$ can be written as
	\begin{equation}\label{25}
		\frac{dL_s}{dz}+\frac{\kappa n_sL_s}{\nu}=\frac{\omega\kappa n_s}{4\pi\nu}G_s(z).
	\end{equation}
	
	The intensity of the basic state can be split into two parts: collimated part $L_s^c$ and diffuse part $L_s^d$ (which occurs due to scattering), such that $L_s=L_s^c+L_s^d$. The equation that governs the collimated part $L_s^c$ is as follows
	\begin{equation}\label{26}
		\frac{dL_s^c}{dz}+\frac{\kappa n_sL_s^c}{\nu}=0,
	\end{equation}
	
	accompanied by boundary conditions

	\begin{equation}\label{27}
		at ~~z=1,~~~~~L_s^c( 1, \theta) =L_t\delta(\boldsymbol{r}-\boldsymbol{r}_0),~~~ where~~~ (\pi/2\leq\theta\leq\pi), 
	\end{equation}
	After calculations, we get 
	\begin{equation}\label{28}
		L_s^c=L_t\exp\left(\int_z^1\frac{\kappa n_s(z')}{\nu}dz'\right)\delta(\boldsymbol{r}-\boldsymbol{r_0}), 
	\end{equation}
	
	and the equation that govern diffused part is   
	\begin{equation}\label{29}
		\frac{dL_s^d}{dz}+\frac{\kappa n_sL_s^d}{\nu}=\frac{\omega\kappa n_s}{4\pi\nu}G_s(z),
	\end{equation}
	accompanied by boundary conditions
	
	\begin{subequations}
		\begin{equation}\label{30a}
			at~~~ z=1,~~~~~~~	L_s^d( z, \theta) =0,~~~ where~~~ (\pi/2\leq\theta\leq\pi), 
		\end{equation}
		\begin{equation}\label{30b}
			at~~~ z=0,~~~~~~~	L_s^d( z, \theta) =0,~~~ where~~~ (0\leq\theta\leq\pi/2). 
		\end{equation}
	\end{subequations}

	Now the total intensity, $G_s=G_s^c+G_s^d$ in the equilibrium state can be written as
	\begin{equation}\label{31}
		G_s^c=\int_0^{4\pi}L_s^c(z,\theta)d\Omega=L_t\exp\left(-\int_z^1\kappa n_s(z')dz'\right),
	\end{equation}
	\begin{equation}\label{32}
		G_s^d=\int_0^{\pi}L_s^d(z,\theta)d\Omega.
	\end{equation}
	
	For no scattering, we get Lambert-Beer law $G_s=G_s^c$. Now, we define a new variable as 
	\begin{equation*}
		\tau=\kappa\int_z^1 n_s(z')dz',
	\end{equation*}
	
	As a result, the total intensity $G_s$ is dependent only on the optical depth $\tau$. In addition, the non-dimensional total intensity, $\Upsilon(\tau)=G_s(\tau)$, can be represented by the following Fredholm Integral Equation (FIE)
	\begin{equation}\label{33}
		\Upsilon(\tau) = e^{-\tau}+\frac{\omega}{2}\int_0^\kappa \Upsilon(\tau')E_1(|\tau-\tau'|)d\tau'.
	\end{equation}
	
	Here, $E_1(x)$ and $E_2(x)$ represent the first and second-order exponential integral, respectively. The FIE has a singularity at $\tau'=\tau$. To solve this FIE, the method of subtraction of singularity is employed.
	
	The radiative flux in the basic state can be expressed as follows
	
	\begin{equation*}
		\boldsymbol{q_s}=\int_0^{4\pi}\left(L_s^c(z,\theta)+L_s^d(z,\theta)\right)\boldsymbol{s}d\Omega=-L_t\exp\left(\int_z^1-\kappa n_s(z')dz'\right)\hat{\boldsymbol{z}}+\int_0^{4\pi}L_s^d(z,\theta)\boldsymbol{s}d\Omega.
	\end{equation*}
	Since $I_s^d(z,\theta)$ is not dependent on $\phi$, the x and y components of $\boldsymbol{q_s}$ become zero. Consequently, we have $\boldsymbol{q}_s=-q_s\hat{\boldsymbol{z}}$, where $q_s=|\boldsymbol{q_s}|$. This implies that the mean swimming direction can be determined as follows
	\begin{equation*}
		<\boldsymbol{p_s}>=-M_s\frac{\boldsymbol{q_s}}{q_s}=M_s\hat{\boldsymbol{z}},
	\end{equation*}

\begin{figure*}[!ht]
	\includegraphics{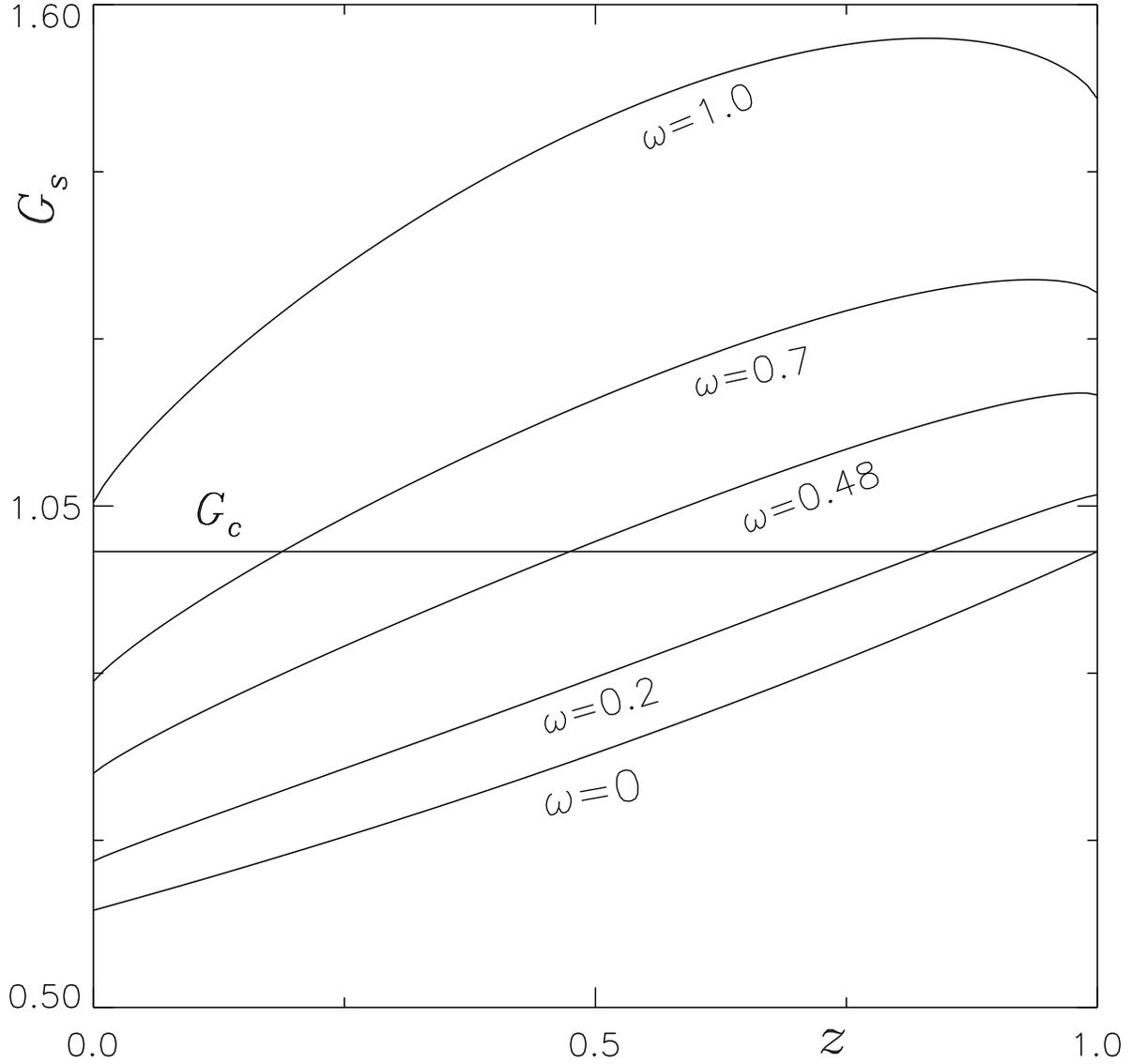}
	\caption{\label{fig2} The total intensity behavior in the uniform suspension is influenced by the fixed parameter values of $S_c=20$, $k=0.5$, and $L_t=1$.}
\end{figure*}

\begin{figure*}[!ht]
	\includegraphics{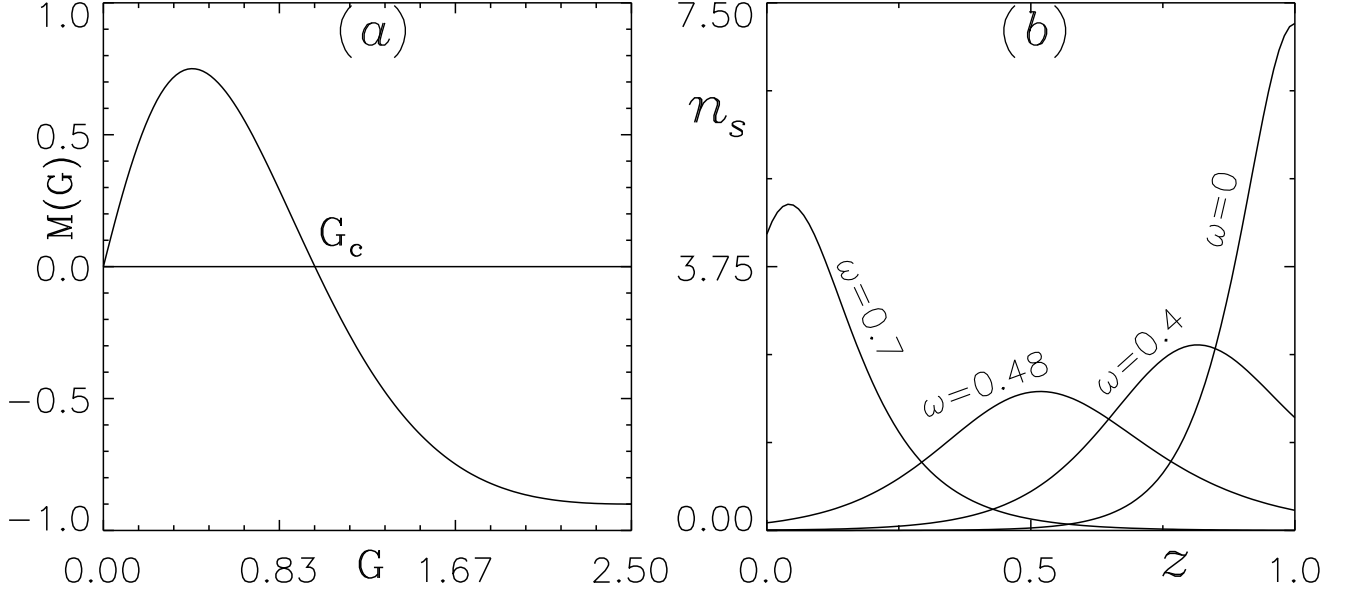}
	\caption{\label{fig3} (a) The taxis responce curve for $G_c=1$ and, (b) the profile at the equilibrium state for different values of scattering albedo. Here, the  governing parameter values $S_c=20,V_c=20,\kappa=0.5$ and $L_t=1$ are kept fixed.}
\end{figure*}
	
	where $M_s=M(G_s).$\par
	The solution for the cell concentration in a basic state can be expressed as $n_s(z)$ and satisfies the following equation
	
	\begin{equation}\label{34}
		\frac{dn_s}{dz}-V_cM_sn_s=0,
	\end{equation}
	where, the basic state cell concentration $n_s(z)$ is accompanied by the conservation relation for the cells
	\begin{equation}\label{35}
		\int_0^1n_s(z)dz=1.
	\end{equation}
	The set of equations given by (\ref{33}) to (\ref{35}) forms a boundary value problem, and the solution to this problem is obtained using a numerical technique called the shooting method.

	In the analysis, we assume that the radiation intensity $L_t=1$ is at the top of the suspension. Here, we also consider the photo response function with $G_c=1$, whose functional form is 
	\begin{equation}\label{36}
		M(G)=0.8\sin\left(\frac{3\pi}{2}\Lambda(G)\right)-0.1\sin\left(\frac{\pi}{2}\Lambda(G)\right),~~~where~~~ \Lambda(G)=\frac{G}{2.5}\exp[0.32(2.5-G)]
	\end{equation}

In Fig.~\ref{fig2}, the variation of the total intensity, $G_s$, is shown across a uniform suspension layer $(n=1)$ for different values of $\omega$ with $\kappa=0.5$. When $\omega=0$, the variation of $G_s$ follows the Lambert-Beer law. However, for non-zero values of $\omega$, the actual $G_s$ at a point is higher than what is predicted by the Lambert-Beer law due to the contribution of scattering.

For $\omega=0$, $\omega=0.2$, and $\omega=0.48$, $G_s$ monotonically decreases across the suspension. In the case of a uniform suspension, the critical intensity, $G_c$, is located at the top ($z=1$) for $\omega=0$ and at around mid-height ($z=0.5$) for $\omega=0.48$. As a result, cells tend to accumulate at the top of the suspension for $\omega=0$ and around the mid-height for $\omega=0.48$ (refer to Fig.~\ref{fig3}(b)).

On the other hand, for $\omega=0.7$, the variation of $G_s$ across the suspension is non-monotonic. In the case of a uniform suspension, the maximum of $I_s$ occurs at $z=0.94$, and the critical intensity, $G_c$, occurs at $z=0.2$. As a result, cells accumulate above the location of $G_c$ if they exhibit negative phototaxis and below it if they exhibit positive phototaxis. In the steady state, the location of the critical intensity ($z=0.1$ for $\omega=0.7$) may not necessarily coincide with that of the uniform state ($z=0.2$ for $\omega=0.7$).

	\section{Linear stability of the problem}
	For stability analysis, we use linear perturbation theory. Here, the small perturbation of amplitude ${\bar{\epsilon}} (0<{\bar{\epsilon}}<<1)$ is made to the equilibrium state, according to the following equation
	\begin{widetext}
		
		\begin{align}\label{37}
			\nonumber[\boldsymbol{u},\zeta,n,L,<p>]=[0,\zeta_s,n_s,L_s,<p_s>]+\bar{\epsilon} [\boldsymbol{u}_1,\zeta_1,n_1,L_1,<\boldsymbol{p}_1>]+\mathcal{O}(\bar{\epsilon}^2)\\
			=[0,\zeta_s,n_s,L_s,<p_s>]+\bar{\epsilon} [\boldsymbol{u}_1,\zeta_1,n_1,L_1^c+L_1^d,<\boldsymbol{p}_1>]+\mathcal{O}(\bar{\epsilon}^2).  
		\end{align}
	\end{widetext}
	Eqs.~(\ref{11})-(\ref{13}) are linearized by substituting the perturbed variables and collecting $o(\bar{\epsilon})$ terms about the equilibrium state, gives
	\begin{equation}\label{38}
		\boldsymbol{\nabla}\cdot \boldsymbol{u}_1=0,
	\end{equation}
	where  $\boldsymbol{u}_1=(u_1,v_1,w_1)$.
	\begin{equation}\label{39}
		Sc^{-1}\left(\frac{\partial \boldsymbol{u_1}}{\partial t}\right)+\sqrt{T_a}(z\times u_1)+\boldsymbol{\nabla} P_{e}+Rn_1\hat{\boldsymbol{z}}=\nabla^{2}\boldsymbol{ u_1},
	\end{equation}
	\begin{equation}\label{40}
		\frac{\partial{n_1}}{\partial{t}}+V_c\boldsymbol{\nabla}\cdot(<\boldsymbol{p_s}>n_1+<\boldsymbol{p_1}>n_s)+w_1\frac{dn_s}{dz}=\boldsymbol{\nabla}^2n_1.
	\end{equation}
	If $G=G_s+\bar{\epsilon} G_1+\mathcal{O}(\bar{\epsilon}^2)=(G_s^c+\bar{\epsilon} G_1^c)+(G_s^d+\bar{\epsilon} G_1^d)+\mathcal{O}(\bar{\epsilon}^2)$, then the collimated total intensity in the basic state is perturbed as $L_t\exp\left(-\kappa\int_z^1(n_s(z')+\bar{\epsilon} n_1+\mathcal{O}(\bar{\epsilon}^2))dz'\right)$  and after simplification, we get
	\begin{equation}\label{41}
		G_1^c=L_t\exp\left(-\int_z^1 \kappa n_s(z')dz'\right)\left(\int_1^z\kappa n_1 dz'\right),
	\end{equation}
	and, $G_1^d$ is given by 
	\begin{equation}\label{42}
		G_1^d=\int_0^{4\pi}L_1^d(\boldsymbol{ x},\boldsymbol{ s})d\Omega.
	\end{equation}
	In the similar manner, the radiative heat flux $q=q_s+\bar{\epsilon} q_1++\mathcal{O}(\bar{\epsilon}^2)=(q_s^c+\bar{\epsilon} q_1^c)+(q_s^d+\bar{\epsilon} q_1^d)+\mathcal{O}(\bar{\epsilon}^2)$, we find 
	\begin{equation}\label{43}
		\boldsymbol{q}_1^c=L_t\exp\left(-\int_z^1 \kappa n_s(z')dz'\right)\left(\int_1^z\kappa n_1 dz'\right)\hat{z}
	\end{equation}
	and
	\begin{equation}\label{44}
		q_1^d=\int_0^{4\pi}L_1^d(\boldsymbol{ x},\boldsymbol{ s})\boldsymbol{ s}d\Omega.
	\end{equation}
	Now the expression 
	\begin{equation*}
		-M(G_s+\bar{\epsilon} G_1)\frac{\boldsymbol{q}_s+\bar{\epsilon}\boldsymbol{q}_1+\mathcal{O}(\bar{\epsilon}^2)}{|\boldsymbol{q}_s+\bar{\epsilon}\boldsymbol{q}_1+\mathcal{O}(\bar{\epsilon}^2)|}-M_s\hat{\boldsymbol{z}},
	\end{equation*}
	provides the perturbed swimming direction when collecting terms of order $O(\bar{\epsilon})$
	\begin{equation}\label{45}
		<\boldsymbol{p_1}>=G_1\frac{dM_s}{dG}\hat{\boldsymbol{z}}-M_s\frac{\boldsymbol{q_1}^H}{\boldsymbol{q_s}},
	\end{equation}
	where $\boldsymbol{q}_1^H=[\boldsymbol{q}_1^x,\boldsymbol{q}_1^y]$ is the perturbed radiative flux $\boldsymbol{q}_1$ in horizontal direction. Substituting the value of $<\boldsymbol{p_1}>$ from Eq.$(\ref{42})$ into Eq.$(\ref{37})$ and simplifying, gives
	\begin{equation}\label{46}
		\frac{\partial{n_1}}{\partial{t}}+V_c\frac{\partial}{\partial z}\left(M_sn_1+n_sG_1\frac{dM_s}{dG}\right)-V_cn_s\frac{M_s}{q_s}\left(\frac{\partial q_1^x}{\partial x}+\frac{\partial q_1^y}{\partial y}\right)+w_1\frac{dn_s}{dz}=\nabla^2n_1.
	\end{equation}
	By elimination of $P_e$ and horizontal component of $u_1$, Eqs. (26), (27) and Eq. (31) can be reduced to three equations for the perturbed variables namely the vertical component of the velocity $w_1$, the vertical component of the vorticity $\zeta_1 (= \zeta\cdot\hat{\boldsymbol{z}})$ and
	the concentration $n_1$. These variables can be decomposed into normal modes as
	\begin{equation}\label{47}
		[w_1,\zeta_1,n_1]=[W(z),Z(z),N(z)]\exp{(\sigma t+i(lx+my))}.  
	\end{equation}
	$L_1^d$ is governed by the equation 
	\begin{equation}\label{48}
		\xi\frac{\partial L_1^d}{\partial x}+\eta\frac{\partial L_1^d}{\partial y}+\nu\frac{\partial L_1^d}{\partial z}+\kappa n_sL_1^d=\frac{\omega\kappa}{4\pi}(n_sG_1^c+n_sG_1^d+G_sn_1)-\kappa n_1L_s,
	\end{equation}
	with the boundary conditions
	\begin{subequations}
		\begin{equation}\label{49a}
			at~~~z=0,1,~~~~~~	L_1^d(x, y, z, \xi, \eta, \nu) =0,~~~where~~~ (\pi/2\leq\theta\leq\pi,~~0\leq\phi\leq 2\pi), 
		\end{equation}
		\begin{equation}\label{49b}
			at~~~z=0,~~~~~~L_1^d(x, y, z,\xi, \eta, \nu) =0,~~~where~~~ (0\leq\theta\leq\pi/2,~~0\leq\phi\leq 2\pi). 
		\end{equation}
	\end{subequations}
	The possible expression for $L_1^d$ suggested by the form of Eq.~$(\ref{45})$ is
	\begin{equation*}
		L_1^d=\Psi_1^d(z,\xi,\eta,\nu)\exp{(\sigma t+i(lx+my))}. 
	\end{equation*}
	From Eqs.~(\ref{38}) and (\ref{39}), we get
	\begin{equation}\label{50}
		G_1^c=\left[L_t\exp\left(-\int_z^1 \kappa n_s(z')dz'\right)\left(\int_1^z\kappa n_1 dz'\right)\right]\exp{(\sigma t+i(lx+my))}=\mathcal{G}_1^c(z)\exp{(\sigma t+i(lx+my))},
	\end{equation}
	and 
	\begin{equation}\label{51}
		G_1^d=\mathcal{G}_1^d(z)\exp{(\sigma t+i(lx+my))}= \left(\int_0^{4\pi}\Psi_1^d(z,\xi,\eta,\nu)d\Omega\right)\exp{(\sigma t+i(lx+my))},
	\end{equation}

	where $\mathcal{G}_1(z)=\mathcal{G}_1^c(z)+\mathcal{G}_1^d(z)$ is the perturbed total intensity.
	
	Now $\Psi_1^d$ satisfies
	\begin{equation}\label{52}
		\frac{d\Psi_1^d}{dz}+\frac{(i(l\nu_1+m\eta)+\kappa n_s)}{\nu}\Psi_1^d=\frac{\omega\kappa}{4\pi\nu}(n_s\mathcal{G}_1+G_s\Theta)-\frac{\kappa}{\nu}L_s\Theta,
	\end{equation}
	with the boundary conditions
	\begin{subequations}
		\begin{equation}\label{53a}
			at~~~z=1,~~~~~~\Psi_1^d( z, \xi, \eta, \nu) =0,~~~where~~~ (\pi/2\leq\theta\leq\pi,~~0\leq\phi\leq 2\pi), 
		\end{equation}
		\begin{equation}\label{53b}
			at~~~z=0,~~~~~~\Psi_1^d( z,\xi, \eta, \nu) =0,~~~where~~~ (0\leq\theta\leq\pi/2,~~0\leq\phi\leq 2\pi). 
		\end{equation}
	\end{subequations}
	In the same manner from Eq.~(\ref{43}), we have
	\begin{equation*}
		q_1^H=[q_1^x,q_1^y]=[P(z),Q(z)]\exp{[\sigma t+i(lx+my)]},
	\end{equation*}
	where
	\begin{equation*}
		P(z)=\int_0^{4\pi}\Psi_1^d(z,\xi,\eta,\nu)\xi d\Omega,\quad Q(z)=\int_0^{4\pi}\Psi_1^d(z,\xi,\eta,\nu)\eta d\Omega.
	\end{equation*}
	The linear stability equations become
	\begin{equation}\label{54}
		\left(\sigma S_c^{-1}+k^2-\frac{d^2}{dz^2}\right)\left( \frac{d^2}{dz^2}-k^2\right)W=Rk^2N,
	\end{equation}
	\begin{equation}\label{55}
	\left(\sigma S_c^{-1}+k^2-\frac{d^2}{dz^2}\right)Z(z)=\sqrt{T_a}\frac{dW}{dz}
	\end{equation}
	\begin{equation}\label{56}
		\left(\sigma+k^2-\frac{d^2}{dz^2}\right)N+V_c\frac{d}{dz}\left(M_sN+n_s\mathcal{G}_1\frac{dM_s}{dG}\right)-i\frac{V_cn_sM_s}{q_s}(lP+mQ)=-\frac{dn_s}{dz}W,
	\end{equation}
	subject to the boundary conditions
	\begin{equation}\label{57}
		at~~~z=0,~~~~~~W=\frac{dW}{dz}=Z(z)=\frac{dN}{dz}-V_cM_sN-n_sV_C\mathcal{G}_1\frac{dM_s}{dG}=0.
	\end{equation}
		\begin{equation}\label{58}
	at~~~z=1,~~~~~~W=\frac{dW}{dz}=\frac{dZ(z)}{dz}=\frac{dN}{dz}-V_cM_sN-n_sV_C\mathcal{G}_1\frac{dM_s}{dG}=0.
	\end{equation}
	
	Here, $k=\sqrt{(l^2+m^2)}$ is the non-dimensional wavenumber.
	
	Eq.~(\ref{52}) becomes (writing D = d/dz)
	\begin{equation}\label{59}
		\aleph_0(z)+\aleph_1(z)\int_1^zN dz+(\sigma+k^2+\aleph_2(z))N+\aleph_3(z)DN-D^2N=-Dn_sW, 
	\end{equation}
	where
	\begin{subequations}
		\begin{equation}\label{60a}
			\aleph_0(z)=V_cD\left(n_s\mathcal{G}_1^d\frac{dM_s}{dG}\right)-i\frac{V_cn_sM_s}{q_s}(lP+mQ),
		\end{equation}
		\begin{equation}\label{60b}
			\aleph_1(z)=\kappa V_cD\left(n_sG_s^c\frac{dM_s}{dG}\right),
		\end{equation}
		\begin{equation}\label{60c}
			\aleph_2(z)=2\kappa V_c n_s G_s^c\frac{dM_s}{dG}+V_c\frac{dM_s}{dG}DG_s^d,
		\end{equation}
		\begin{equation}\label{60d}
			\aleph_3(z)=V_cM_s.
		\end{equation}
	\end{subequations}
	Now, introducing a new variable as
	\begin{equation}\label{61}
		\Phi(z)=\int_1^zN(z')dz',
	\end{equation}
	Eq.~\ref{51} and Eq.~\ref{54} becoems
	\begin{equation}\label{62}
		\left(\sigma S_c^{-1}+k^2-D^2\right)\left( D^2-k^2\right)W=Rk^2D\Phi,
	\end{equation}
		\begin{equation}\label{63}
	\left(\sigma S_c^{-1}+k^2-D^2\right)Z(z)=\sqrt{T_a}DW
	\end{equation}
	\begin{equation}\label{64}
		\aleph_0(z)+\aleph_1(z)\Phi+(\sigma+k^2+\aleph_2(z))D\Phi+\aleph_3(z)D^2\Phi-D^3\Phi=-Dn_sW, 
	\end{equation}
	with the boundary conditions,
	
	\begin{equation}\label{65}
	at~~~z=0,~~~~~~W=DW=Z(z)=D^2\Phi-V_cM_s\Phi-n_sV_c\mathcal{G}_1\frac{dM_s}{dG}=0.
	\end{equation}
	\begin{equation}\label{66}
	at~~~z=1,~~~~~~W=DW=DZ(z)=D^2\Phi-V_cM_s\Phi-n_sV_C\mathcal{G}_1\frac{dM_s}{dG}=0,
	\end{equation}

	and the additional boundary condition is,
	\begin{equation}\label{67}
		at~~~z=1,~~~~~~\Phi(z)=0.
	\end{equation}

	\section{SOLUTION technique}
	In order to find solutions for Eqs. (\ref{62}) to (\ref{64}), we utilize the Newton-Raphson-Kantorovich (NRK) iterative method, as described in the work by Cash et al.~\cite{19cash1980}. This numerical method enables us to calculate the growth rate, Re$(\sigma)$, or neutral stability curves in the $(k, R)$-plane for a specific parameter set.	
	The neutral curve, denoted as $R^{(n)}(k)$, where $n$ is an integer greater than or equal to 1, consists of an infinite number of branches. It provides a unique solution to the linear stability problem for a given set of parameters. Among these branches, the one with the lowest value of $R$ is considered the most significant, and the corresponding bioconvective solution is identified as $(k_c, R_c)$. This particular solution is referred to as the most unstable solution.	
	By utilizing the equation $\lambda_c=2\pi/k_c$, where $\lambda_c$ represents the wavelength of the initial disturbance, we can determine the wavelength associated with the most unstable solution.
	
	\section{NUMERICAL RESULTS}
	\begin{figure*}[!ht]
		\includegraphics{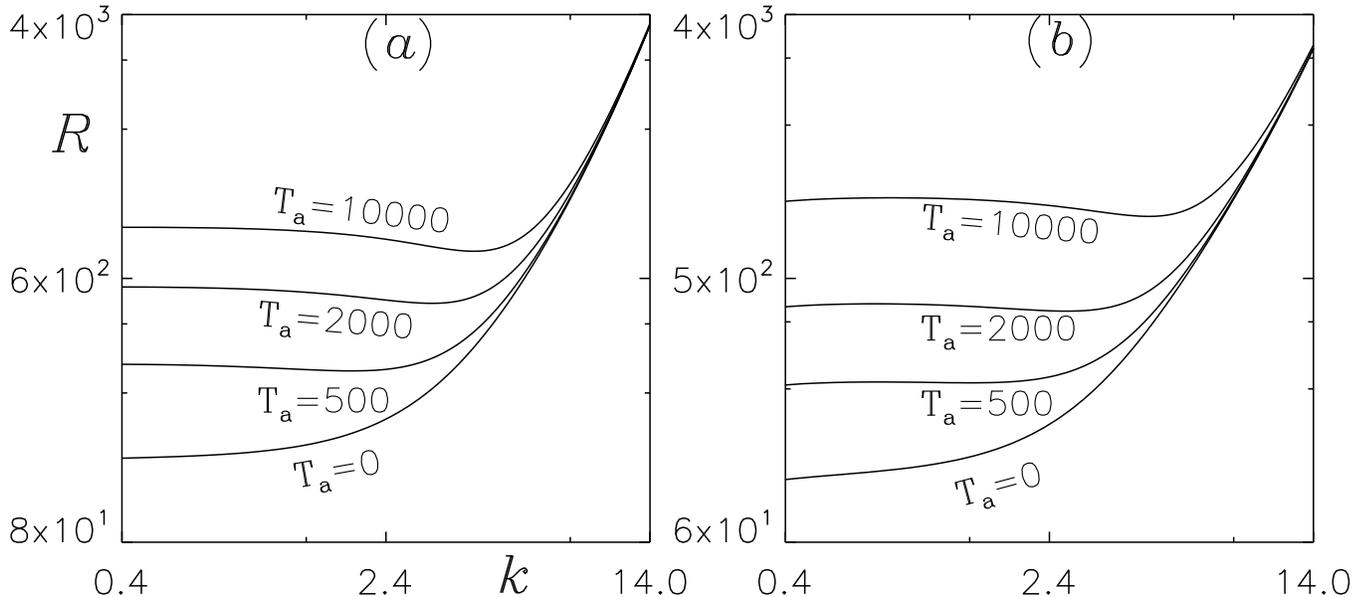}
		\caption{\label{fig4} The marginal stability curves (a) for $\kappa=0.5$ and, (b) for $\omega=1$. Here, the parameter values $S_c=20,V_c=20 L_t=1$, and $\omega=0$ are kept fixed.}
	\end{figure*}
	
In this study, we conduct a systematic analysis to examine the effect of rotation, specifically the Taylor number $T_a$, while keeping the parameters $S_c$, $L_t$, $V_c$, $\kappa$, and $\omega$ constant. With a wide range of parameter values to consider, comprehensively exploring the entire parameter space can be challenging. Therefore, we focused on investigating the impact of a discrete set of constant parameter values on the onset of bioconvection.
Throughout the study, we maintain $S_c=20$ and $L_t=1$ as fixed parameter values. The parameters related to scattering albedo, extinction coefficient, and cell swimming speed were set as follows: $\omega$ ranged from 0 to 1, $\kappa$ took values of 0.5 and 1.0, and $V_c$ varied among 10, 15, and 20. These constant parameter values is chosen to examine their specific influence on the initiation of bioconvection in the system under investigation.

\begin{figure*}[!ht]
	\includegraphics{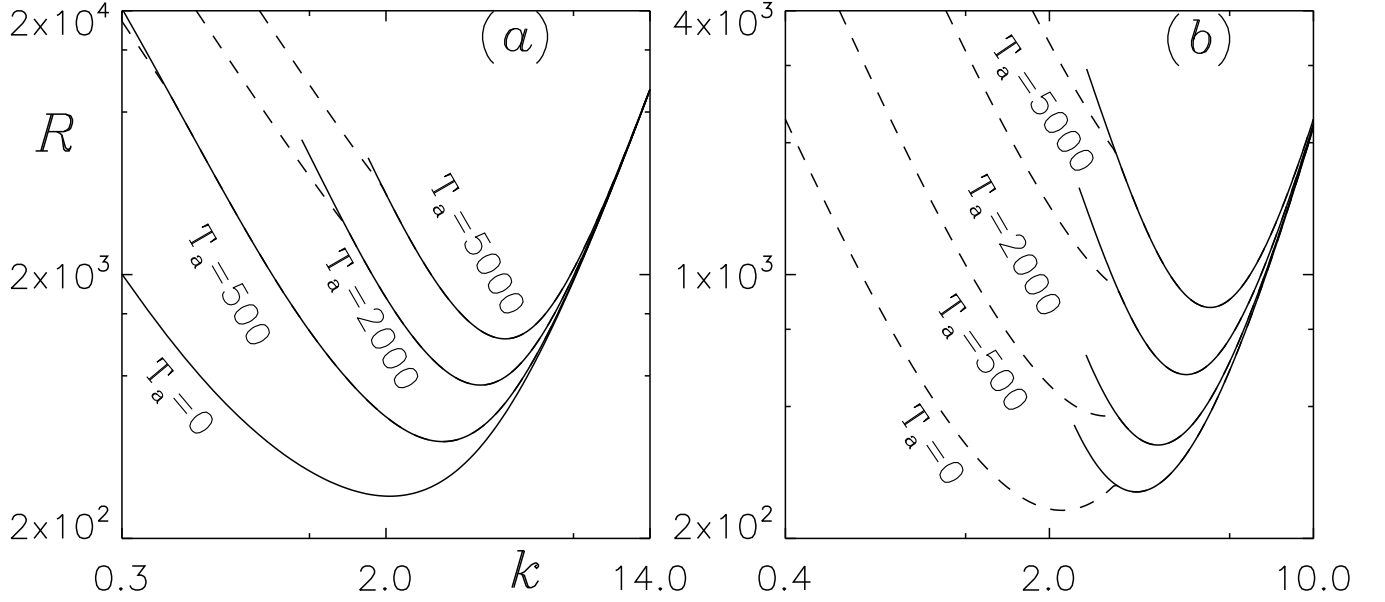}
	\caption{\label{fig5}The marginal stability curves (a) for $\kappa=0.5, \omega=0.43$ and, (b) for $\omega=1,\omega=0.59$. Here, the parameter values $S_c=20,V_c=20$, and $L_t=1$ are kept fixed.}
\end{figure*}
\begin{figure*}[!ht]
	\includegraphics{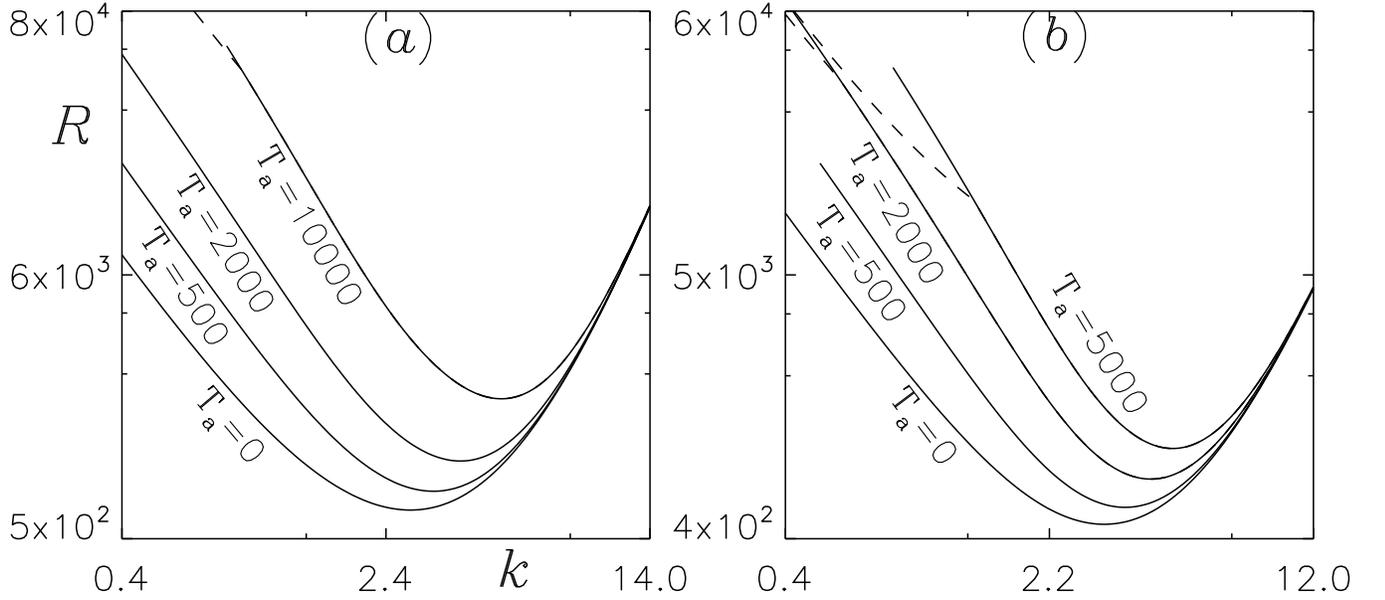}
	\caption{\label{fig6} The marginal stability curves (a) for $\kappa=0.5,\omega=0.48$ and, (b) for $\omega=1,\omega=0.61$. Here, the parameter values $S_c=20,V_c=20$, and $L_t=1$ are kept fixed.}
\end{figure*}
	
	\subsection{ WHEN SCATTERING IS WEAK}
	In this study, our focus is on investigating the impact of rotation on the initiation of bioconvection under the condition of weak scattering. To ensure strong self-shading effects, we specifically selected a lower value for the scattering albedo, denoted as $\omega$. Two different scenarios were considered based on the value of the extinction coefficient, $\kappa$. When $\kappa$ is equal to 0.5, the self-shading is categorized as weak, whereas a value of 1 indicates strong self-shading.
	
	Throughout this section, we consistently employed a critical intensity of $G_c=1$. This value was utilized as a reference point to analyze the effects of rotation on the initiation of bioconvection in the presence of weak scattering and varying degrees of self-shading.

\subsubsection{$V_c=20$}
	
In this section, we investigate the influence of the Taylor number, denoted as $T_a$, on the bioconvective instability at three distinct locations within the suspension's equilibrium state. Specifically, we examine the top ($z\approx 1$), three-quarter ($z\approx 3/4$), and mid-height ($z\approx 1/2$) positions of the suspension. The parameters $V_c=20$ and $\kappa=0.5$ (Case 1) and $\kappa=1$ (Case 2) are held fixed throughout the analysis.

We systematically vary the value of $T_a$ from 0 to higher nonzero values, up to 10000, in order to observe the resulting changes in the bioconvective instability. By exploring different values of $T_a$, we aim to understand how rotation affects the onset and characteristics of bioconvection at different locations within the suspension. The two cases considered, differing in the value of $\kappa$, allow us to analyze the influence of the extinction coefficient on the bioconvective behavior under the influence of rotation.

\begin{table}[!htbp]
	\caption{\label{tab1}The table shows the numerical results of bioconvective solutions for a constant $V_c=$20 and $\kappa=$0.5. The values are presented for different increments in $T_a$, while keeping all other parameters unchanged. The presence of a star symbol indicates that the $R^{(1)}(k)$ branch of the neutral curve exhibits oscillatory behavior.}
	\begin{ruledtabular}
		\begin{tabular}{cccccccccc}
			\multirow{2}{*}{$T_a$}&\multicolumn{3}{c}{$\omega=0$} & \multicolumn{3}{c}{$\omega=0.43$} & \multicolumn{3}{c}{$\omega=0.48$}\\\cline{2-4}\cline{5-7}\cline{8-10} & $\lambda_c$ & $R_c$ & $Im(\sigma)$ &$ \lambda_c$ & $R_c$ & $Im(\sigma)$ &$\lambda_c$ & $R_c$ & $Im(\sigma)$\\
			\hline
			0 & $\infty$ & 139.83 & 0 & 2.98 & 221.85 & 0 & 2.26 & 658.78 & 0 \\ 		
			100 & $\infty$ & 177.98 & 0 & 2.63 & 259.30 & 0 & 2.15 & 689.10 & 0 \\ 
			500 & 3.34 & 284.96 & 0 & 2.01$^{\star}$ & 368.14 & 0 & 1.94 & 789.35 & 0 \\ 
			1000 & 2.51 & 361.02 & 0 & 1.77$^{\star}$ & 467.87 & 0 & 1.78 & 890.33 & 0 \\ 
			2000 & 1.94 & 470.27 & 0 & 1.54$^{\star}$ & 621.83 & 0 & 1.61 & 1035.63 & 0 \\ 
			5000 & 1.48 & 691.95 & 0 & 1.29$^{\star}$ & 954.64 & 0 & 1.39 & 1434.92 & 0 \\ 
			10000 & 1.22 & 952.34 & 0 & 1.14$^{\star}$ & 1361.65 & 0 & 1.22$^{\star}$ & 1921.36 & 0 \\ 
			
		\end{tabular}
	\end{ruledtabular}
	
\end{table}

\begin{table}[!htbp]
	\caption{\label{tab2}The table shows the numerical results of bioconvective solutions for a constant $V_c=$20 and $\kappa=$1. The values are presented for different increments in $T_a$, while keeping all other parameters unchanged. The presence of a star symbol indicates that the $R^{(1)}(k)$ branch of the neutral curve exhibits oscillatory behavior, while the double dagger symbol indicates the presence of an overstable solution.}
	\begin{ruledtabular}
		\begin{tabular}{cccccccccc}
			\multirow{2}{*}{$T_a$}&\multicolumn{3}{c}{$\omega=0$} & \multicolumn{3}{c}{$\omega=0.59$} & \multicolumn{3}{c}{$\omega=0.61$}\\\cline{2-4}\cline{5-7}\cline{8-10} & $\lambda_c$ & $R_c$ & $Im(\sigma)$ &$ \lambda_c$ & $R_c$ & $Im(\sigma)$ &$\lambda_c$ & $R_c$ & $Im(\sigma)$\\
			\hline
			0 & $\infty$ & 25.05 & 0 & 2.90$^{\ddag}$ & 289.72$^{\ddag}$ & 8.06 & 2.01 & 457.79 & 0 \\ 		
			100 & $\infty$ & 32.46 & 0 & 2.69$^{\ddag}$ & 332.96$^{\ddag}$ & 7.92 & 1.94 & 475.81 & 0 \\ 
			500 & 4.51 & 213.39 & 0 & 1.63$^{\star}$ & 408.41 & 0 & 1.76 & 536.22 & 0 \\ 
			1000 & 3.12 & 280.32 & 0 & 1.50$^{\star}$ & 591.26 & 0 & 1.64 & 597.51 & 0 \\ 
			2000 & 2.35 & 377.55 & 0 & 1.37$^{\star}$ & 841.32 & 0 & 1.50$^{\star}$ & 698.06 & 0 \\ 
			5000 & 1.71 & 574.57 & 0 & 1.78$^{\star}$ & 954.64 & 0 & 1.29$^{\star}$ & 928.46 & 0 \\ 
			10000 & 1.37 & 802.62 & 0 & 1.05$^{\star}$ & 1146.14 & 0 & 1.16$^{\star}$ & 1222.62 & 0 \\ 
			
		\end{tabular}
	\end{ruledtabular}
	
\end{table}	
	
In Fig.\ref{fig4}(a), the marginal stability curves for $\kappa=0.5$ are presented, while Fig.\ref{fig4}(b) shows the marginal stability curves for $\kappa=1$. These curves correspond to different levels of the Taylor number ($T_a$), while maintaining constant values of $V_c=20$ and $\omega=0$.
When $T_a=0$, the linear stability analysis predicts a stationary solution at the bioconvective instability. In this case, the critical wavelength is infinite. As we increase the Taylor number to 100, the solution remains stationary with an infinite critical wavelength, and the critical Rayleigh number also increases. When $T_a=500$, the critical wavelength becomes finite, and the critical Rayleigh number continues to increase. As the Taylor number further increases up to 10000, the critical wavelength decreases, and the critical Rayleigh number increases.

In Fig.\ref{fig5}(a), the marginal stability curve for $\kappa=0.5$ and $\omega=0.43$ is shown, while Fig.\ref{fig5}(b) depicts the marginal stability curve for $\kappa=1$ and $\omega=0.59$. In both cases, the equilibrium state is formed at the three-quarter height of the domain.

For $T_a=0$, the marginal stability curve exhibits a finite critical wavelength and a finite critical Rayleigh number. In the case of $\kappa=0.5$, a stationary solution is observed, whereas for $\kappa=1$, an oscillatory branch bifurcates from the stationary branch of the marginal stability curve, indicating the presence of an overstable solution. As the Taylor number increases to 100, the marginal stability curve displays a behavior similar to that observed at $T_a=0$. When $T_a=500$, for $\kappa=0.5$, an oscillatory branch bifurcates from the stationary branch of the marginal stability curve, but the critical solution occurs at the stationary branch of the neutral curve. The same behavior is observed for $\kappa=1$, where the overstable solution converts into a stationary solution. For higher values of the Taylor number $T_a$, the same behavior is observed. In both cases, the critical wavelength decreases, and the critical Rayleigh number increases as the Taylor number increases.

\begin{table}[!htbp]
	\caption{\label{tab3}The table shows the numerical results of bioconvective solutions for a constant $V_c=$10 and $\kappa=$0.5. The values are presented for different increments in $T_a$, while keeping all other parameters unchanged.}
	\begin{ruledtabular}
		\begin{tabular}{cccccccccc}
			\multirow{2}{*}{$T_a$}&\multicolumn{3}{c}{$\omega=0$} & \multicolumn{3}{c}{$\omega=0.37$} & \multicolumn{3}{c}{$\omega=0.47$}\\\cline{2-4}\cline{5-7}\cline{8-10} & $\lambda_c$ & $R_c$ & $Im(\sigma)$ &$ \lambda_c$ & $R_c$ & $Im(\sigma)$ &$\lambda_c$ & $R_c$ & $Im(\sigma)$\\
			\hline
			0 & $\infty$ & 103.46 & 0 & 7.76 & 210.42 & 0 & 2.22 & 1757.31 & 0 \\ 		
			100 & 10.76 & 140.29 & 0 & 4.98 & 284.56 & 0 & 2.13 & 1816.54 & 0 \\ 
			500 & 3.15 & 237.77 & 0 & 2.63 & 485.66 & 0 & 1.97 & 2022.56 & 0 \\ 
			1000 & 2.40 & 324.96 & 0 & 2.13 & 656.80 & 0 & 1.81 & 2238.98 & 0 \\ 
			2000 & 1.91 & 461.62 & 0 & 1.76 & 912.87 & 0 & 1.66 & 2601.61 & 0 \\ 
			5000 & 1.46 & 765.30 & 0 & 1.42 & 1458.40 & 0 & 1.44 & 3447.30 & 0 \\ 
			10000 & 1.22 & 1146.71 & 0 & 1.23 & 2127.16 & 0 & 1.26 & 4534.31 & 0 \\ 
			
		\end{tabular}
	\end{ruledtabular}
	
\end{table}

\begin{table}[!htbp]
	\caption{\label{tab4}The table shows the numerical results of bioconvective solutions for a constant $V_c=$20 and $\kappa=$1. The values are presented for different increments in $T_a$, while keeping all other parameters unchanged. The presence of a star symbol indicates that the $R^{(1)}(k)$ branch of the neutral curve exhibits oscillatory behavior.}
	\begin{ruledtabular}
		\begin{tabular}{cccccccccc} 
			\multirow{2}{*}{$T_a$}&\multicolumn{3}{c}{$\omega=0$} &
			\multicolumn{3}{c}{$\omega=0.53$} &
			\multicolumn{3}{c}{$\omega=0.6$}\\\cline{2-4}\cline{5-7}\cline{8-10} & $\lambda_c$ & $R_c$ & $Im(\sigma)$ &$ \lambda_c$ & $R_c$ & $Im(\sigma)$ &$\lambda_c$ & $R_c$ & $Im(\sigma)$\\
			\hline
			0 & $\infty$ & 82.27 & 0 & 3.44 & 215.32 & 0 & 2.27 & 778.29 & 0 \\ 		
			100 & $\infty$ & 111.20 & 0 & 2.85 & 256.98 & 0 & 2.14 & 814.20 & 0 \\ 
			500 & 3.67 & 196.55 & 0 & 2.14 & 376.86 & 0 & 1.94 & 932.72 & 0 \\ 
			1000 & 2.69 & 270.12 & 0 & 1.86 & 486.79 & 0 & 1.79 & 1052.00 & 0 \\ 
			2000 & 2.05 & 384.05 & 0 & 1.59$^{\star}$ & 656.98 & 0 & 1.61 & 1247.15 & 0 \\ 
			5000 & 1.54 & 632.28 & 0 & 1.32$^{\star}$ & 1026.73 & 0 & 1.39 & 1695.11 & 0 \\ 
			10000 & 1.26 & 938.37 & 0 & 1.16$^{\star}$ & 1481.61 & 0 & 1.23 & 2270.06 & 0 \\ 
			
		\end{tabular}
	\end{ruledtabular}
	
\end{table}

Fig.\ref{fig6}(a) and Fig.\ref{fig6}(b) illustrate the marginal stability curves for $\kappa=0.5$, $\omega=0.48$, and $\kappa=1$, $\omega=61$, respectively. In these cases, the accumulation of cells occurs at the mid-height of the domain.

When $T_a=0$, both cases ($\kappa=0.5$ and $\kappa=1$) exhibit finite values for the critical wavelength and critical Rayleigh number. For $\kappa=0.5$, when the Taylor number increases to 5000, the same behavior is observed, but at $T_a=10000$, an oscillatory branch bifurcates from the marginal stability curve. However, the most unstable solution still occurs at the stationary branch, resulting in a stationary mode of bio-convective instability. In the case of $\kappa=1$, an oscillatory branch bifurcates from the stationary branch of the marginal stability curve for $T_a=2000$, but the most unstable solution occurs at the stationary branch of the neutral curve, resulting in a stationary solution. When $T_a$ is further increased to 10000, the same behavior is observed. In both cases, the critical wavelength decreases, and the critical Rayleigh number increases as the Taylor number increases.

\begin{table}[!htbp]
	\caption{\label{tab5}The table shows the numerical results of bioconvective solutions for a constant $V_c=$15 and $\kappa=$0.5. The values are presented for different increments in $T_a$, while keeping all other parameters unchanged. The presence of a star symbol indicates that the $R^{(1)}(k)$ branch of the neutral curve exhibits oscillatory behavior.}
	\begin{ruledtabular}
		\begin{tabular}{cccccccccc}
			\multirow{2}{*}{$T_a$}&\multicolumn{3}{c}{$\omega=0$} & \multicolumn{3}{c}{$\omega=0.4$} & \multicolumn{3}{c}{$\omega=0.475$}\\\cline{2-4}\cline{5-7}\cline{8-10} & $\lambda_c$ & $R_c$ & $Im(\sigma)$ &$ \lambda_c$ & $R_c$ & $Im(\sigma)$ &$\lambda_c$ & $R_c$ & $Im(\sigma)$\\
			\hline
			0 & $\infty$ & 120.43 & 0 & 4.39 & 191.08 & 0 & 2.23 & 979.51 & 0 \\ 		
			100 & $\infty$ & 157.24 & 0 & 3.44 & 240.03 & 0 & 2.14 & 1019.07 & 0 \\ 
			500 & 3.15 & 249.30 & 0 & 2.30 & 375.20 & 0 & 1.94 & 1152.62 & 0 \\ 
			1000 & 3.15 & 331.43 & 0 & 1.93 & 495.25 & 0 & 1.81 & 1289.57 & 0 \\ 
			2000 & 1.91 & 439.63 & 0 & 1.63 & 678.40 & 0 & 1.65 & 1516.15 & 0 \\ 
			5000 & 1.46 & 682.20 & 0 & 1.35$^{\star}$ & 1072.65 & 0 & 1.42 & 2039.68 & 0 \\ 
			10000 & 1.20 & 977.53 & 0 & 1.18$^{\star}$ & 1556.43 & 0 & 1.24 & 2712.74 & 0 \\ 
			
		\end{tabular}
	\end{ruledtabular}
	
\end{table}

\begin{table}[!htbp]
	\caption{\label{tab6}The table shows the numerical results of bioconvective solutions for a constant $V_c=$20 and $\kappa=$1. The values are presented for different increments in $T_a$, while keeping all other parameters unchanged. The presence of a star symbol indicates that the $R^{(1)}(k)$ branch of the neutral curve exhibits oscillatory behavior.}
	\begin{ruledtabular}
		\begin{tabular}{cccccccccc} 
			\multirow{2}{*}{$T_a$}&\multicolumn{3}{c}{$\omega=0$} &
			\multicolumn{3}{c}{$\omega=0.57$} &
			\multicolumn{3}{c}{$\omega=0.61$}\\\cline{2-4}\cline{5-7}\cline{8-10} & $\lambda_c$ & $R_c$ & $Im(\sigma)$ &$ \lambda_c$ & $R_c$ & $Im(\sigma)$ &$\lambda_c$ & $R_c$ & $Im(\sigma)$\\
			\hline
			0 & $\infty$ & 72.32 & 0 & 2.25$^{\star}$ & 267.48 & 0 & 2.13 & 559.66 & 0 \\ 		
			100 & $\infty$ & 95.36 & 0 & 2.08$^{\star}$ & 294.69 & 0 & 2.04 & 582.84 & 0 \\ 
			500 & 4.07 & 198.20 & 0 & 1.82$^{\star}$ & 381.19 & 0 & 1.84 & 660.36 & 0 \\ 
			1000 & 2.90 & 265.33 & 0 & 1.64$^{\star}$ & 465.12 & 0 & 1.72 & 739.25 & 0 \\ 
			2000 & 2.21 & 365.74 & 0 & 1.46$^{\star}$ & 598.11 & 0 & 1.56 & 868.82 & 0 \\ 
			5000 & 1.60 & 575.90 & 0 & 1.24$^{\star}$ & 889.67 & 0 & 1.35$^{\star}$ & 1167.19 & 0 \\ 
			10000 & 1.30 & 826.26 & 0 & 1.10$^{\star}$ & 1246.45 & 0 & 1.20$^{\star}$ & 1550.34 & 0 \\ 
			
		\end{tabular}
	\end{ruledtabular}
	
\end{table}

\subsubsection{$V_c=10$ and $V_c=15$}
	We have also investigated the effect of rotation (i.e. Taylor number) on bio-convective instability for $V_c=10$ and $V_c=15$. The critical Rayleigh number, $R_c$, and wavelength, $\lambda_c$, obtained from the numerical simulations are shown in Table~\ref{tab2}.

	\section{Conclusion}
In this study, we have developed a novel model that considers the influence of rotation on the bio-convective instability in an isotropic scattering suspension of phototactic microorganisms. This is the first time such effects have been incorporated into the model. The linear stability of the suspension has been examined using the linear perturbation theory.

Our findings reveal that the variation of total intensity across the depth of the suspension is non-monotonic due to scattering, resulting in the occurrence of critical intensity at one or two locations within the suspension. Consequently, the microorganisms tend to accumulate around a single horizontal layer or two horizontal layers, depending on the nature of the suspension. Two-layer accumulation is observed in predominantly scattering suspensions.

The linear stability analysis of the suspension predicts both stationary and oscillatory modes of instability at the onset of bioconvection. However, the overstable solution is only observed when there is a combination of higher cell swimming speed, higher extinction coefficient, and low rotation rate (indicated by a low Taylor number). At $\omega=0$, the wavelength of the initial pattern is always infinite at low Taylor numbers. As the Taylor number increases, the critical wavelength decreases, and the critical Rayleigh number increases, indicating enhanced stability of the suspension at higher rotation rates.

It is important to note that while our model provides valuable insights, its validation requires comparison with experimental data on purely phototactic bioconvection. Unfortunately, such data is currently unavailable due to the scarcity of microorganisms exhibiting pure phototactic behavior. Future studies should focus on identifying suitable microorganisms for this purpose. Furthermore, our proposed model can be extended to simulate other phototactic bioconvection phenomena of interest.

	\section*{ Availability of Data}
	The supporting data of this article is available within the article. 
	\nocite{*}
	%	\section*{REFERENCES}
	\bibliography{isotropic_rotation}

\providecommand{\noopsort}[1]{}\providecommand{\singleletter}[1]{#1}%
\begin{thebibliography}{41}%
\makeatletter
\providecommand \@ifxundefined [1]{%
 \@ifx{#1\undefined}
}%
\providecommand \@ifnum [1]{%
 \ifnum #1\expandafter \@firstoftwo
 \else \expandafter \@secondoftwo
 \fi
}%
\providecommand \@ifx [1]{%
 \ifx #1\expandafter \@firstoftwo
 \else \expandafter \@secondoftwo
 \fi
}%
\providecommand \natexlab [1]{#1}%
\providecommand \enquote  [1]{``#1''}%
\providecommand \bibnamefont  [1]{#1}%
\providecommand \bibfnamefont [1]{#1}%
\providecommand \citenamefont [1]{#1}%
\providecommand \href@noop [0]{\@secondoftwo}%
\providecommand \href [0]{\begingroup \@sanitize@url \@href}%
\providecommand \@href[1]{\@@startlink{#1}\@@href}%
\providecommand \@@href[1]{\endgroup#1\@@endlink}%
\providecommand \@sanitize@url [0]{\catcode `\\12\catcode `\$12\catcode
  `\&12\catcode `\#12\catcode `\^12\catcode `\_12\catcode `\%12\relax}%
\providecommand \@@startlink[1]{}%
\providecommand \@@endlink[0]{}%
\providecommand \url  [0]{\begingroup\@sanitize@url \@url }%
\providecommand \@url [1]{\endgroup\@href {#1}{\urlprefix }}%
\providecommand \urlprefix  [0]{URL }%
\providecommand \Eprint [0]{\href }%
\providecommand \doibase [0]{http://dx.doi.org/}%
\providecommand \selectlanguage [0]{\@gobble}%
\providecommand \bibinfo  [0]{\@secondoftwo}%
\providecommand \bibfield  [0]{\@secondoftwo}%
\providecommand \translation [1]{[#1]}%
\providecommand \BibitemOpen [0]{}%
\providecommand \bibitemStop [0]{}%
\providecommand \bibitemNoStop [0]{.\EOS\space}%
\providecommand \EOS [0]{\spacefactor3000\relax}%
\providecommand \BibitemShut  [1]{\csname bibitem#1\endcsname}%
\let\auto@bib@innerbib\@empty
%</preamble>
\bibitem [{\citenamefont {Platt}(1961)}]{20platt1961}%
  \BibitemOpen
  \bibfield  {author} {\bibinfo {author} {\bibfnamefont {J.~R.}\ \bibnamefont
  {Platt}},\ }\bibfield  {title} {\enquote {\bibinfo {title} {" bioconvection
  patterns" in cultures of free-swimming organisms},}\ }\href@noop {}
  {\bibfield  {journal} {\bibinfo  {journal} {Science}\ }\textbf {\bibinfo
  {volume} {133}},\ \bibinfo {pages} {1766--1767} (\bibinfo {year}
  {1961})}\BibitemShut {NoStop}%
\bibitem [{\citenamefont {Pedley}\ and\ \citenamefont
  {Kessler}(1992)}]{21pedley1992}%
  \BibitemOpen
  \bibfield  {author} {\bibinfo {author} {\bibfnamefont {T.~J.}\ \bibnamefont
  {Pedley}}\ and\ \bibinfo {author} {\bibfnamefont {J.~O.}\ \bibnamefont
  {Kessler}},\ }\bibfield  {title} {\enquote {\bibinfo {title} {Hydrodynamic
  phenomena in suspensions of swimming microorganisms},}\ }\href@noop {}
  {\bibfield  {journal} {\bibinfo  {journal} {Annual Review of Fluid
  Mechanics}\ }\textbf {\bibinfo {volume} {24}},\ \bibinfo {pages} {313--358}
  (\bibinfo {year} {1992})}\BibitemShut {NoStop}%
\bibitem [{\citenamefont {Hill}\ and\ \citenamefont
  {Pedley}(2005)}]{22hill2005}%
  \BibitemOpen
  \bibfield  {author} {\bibinfo {author} {\bibfnamefont {N.~A.}\ \bibnamefont
  {Hill}}\ and\ \bibinfo {author} {\bibfnamefont {T.~J.}\ \bibnamefont
  {Pedley}},\ }\bibfield  {title} {\enquote {\bibinfo {title}
  {Bioconvection},}\ }\href@noop {} {\bibfield  {journal} {\bibinfo  {journal}
  {Fluid Dynamics Research}\ }\textbf {\bibinfo {volume} {37}},\ \bibinfo
  {pages} {1} (\bibinfo {year} {2005})}\BibitemShut {NoStop}%
\bibitem [{\citenamefont {Bees}(2020)}]{23bees2020}%
  \BibitemOpen
  \bibfield  {author} {\bibinfo {author} {\bibfnamefont {M.~A.}\ \bibnamefont
  {Bees}},\ }\bibfield  {title} {\enquote {\bibinfo {title} {Advances in
  bioconvection},}\ }\href@noop {} {\bibfield  {journal} {\bibinfo  {journal}
  {Annual Review of Fluid Mechanics}\ }\textbf {\bibinfo {volume} {52}},\
  \bibinfo {pages} {449--476} (\bibinfo {year} {2020})}\BibitemShut {NoStop}%
\bibitem [{\citenamefont {Javadi}\ \emph {et~al.}(2020)\citenamefont {Javadi},
  \citenamefont {Arrieta}, \citenamefont {Tuval},\ and\ \citenamefont
  {Polin}}]{24javadi2020}%
  \BibitemOpen
  \bibfield  {author} {\bibinfo {author} {\bibfnamefont {A.}~\bibnamefont
  {Javadi}}, \bibinfo {author} {\bibfnamefont {J.}~\bibnamefont {Arrieta}},
  \bibinfo {author} {\bibfnamefont {I.}~\bibnamefont {Tuval}}, \ and\ \bibinfo
  {author} {\bibfnamefont {M.}~\bibnamefont {Polin}},\ }\bibfield  {title}
  {\enquote {\bibinfo {title} {Photo-bioconvection: towards light control of
  flows in active suspensions},}\ }\href@noop {} {\bibfield  {journal}
  {\bibinfo  {journal} {Philosophical Transactions of the Royal Society A}\
  }\textbf {\bibinfo {volume} {378}},\ \bibinfo {pages} {20190523} (\bibinfo
  {year} {2020})}\BibitemShut {NoStop}%
\bibitem [{\citenamefont {Wager}(1911)}]{1wager1911}%
  \BibitemOpen
  \bibfield  {author} {\bibinfo {author} {\bibfnamefont {H.}~\bibnamefont
  {Wager}},\ }\bibfield  {title} {\enquote {\bibinfo {title} {Vii. on the
  effect of gravity upon the movements and aggregation of euglena viridis,
  ehrb., and other micro-organisms},}\ }\href@noop {} {\bibfield  {journal}
  {\bibinfo  {journal} {Philosophical Transactions of the Royal Society of
  London. Series B, Containing Papers of a Biological Character}\ }\textbf
  {\bibinfo {volume} {201}},\ \bibinfo {pages} {333--390} (\bibinfo {year}
  {1911})}\BibitemShut {NoStop}%
\bibitem [{\citenamefont {Kitsunezaki}, \citenamefont {Komori},\ and\
  \citenamefont {Harumoto}(2007)}]{2kitsunezaki2007}%
  \BibitemOpen
  \bibfield  {author} {\bibinfo {author} {\bibfnamefont {S.}~\bibnamefont
  {Kitsunezaki}}, \bibinfo {author} {\bibfnamefont {R.}~\bibnamefont {Komori}},
  \ and\ \bibinfo {author} {\bibfnamefont {T.}~\bibnamefont {Harumoto}},\
  }\bibfield  {title} {\enquote {\bibinfo {title} {Bioconvection and front
  formation of paramecium tetraurelia},}\ }\href@noop {} {\bibfield  {journal}
  {\bibinfo  {journal} {Physical Review E}\ }\textbf {\bibinfo {volume} {76}},\
  \bibinfo {pages} {046301} (\bibinfo {year} {2007})}\BibitemShut {NoStop}%
\bibitem [{\citenamefont {Kessler}(1985)}]{3kessler1985}%
  \BibitemOpen
  \bibfield  {author} {\bibinfo {author} {\bibfnamefont {J.~O.}\ \bibnamefont
  {Kessler}},\ }\bibfield  {title} {\enquote {\bibinfo {title} {Co-operative
  and concentrative phenomena of swimming micro-organisms},}\ }\href@noop {}
  {\bibfield  {journal} {\bibinfo  {journal} {Contemporary Physics}\ }\textbf
  {\bibinfo {volume} {26}},\ \bibinfo {pages} {147--166} (\bibinfo {year}
  {1985})}\BibitemShut {NoStop}%
\bibitem [{\citenamefont {Williams}\ and\ \citenamefont
  {Bees}(2011)}]{4williams2011}%
  \BibitemOpen
  \bibfield  {author} {\bibinfo {author} {\bibfnamefont {C.~R.}\ \bibnamefont
  {Williams}}\ and\ \bibinfo {author} {\bibfnamefont {M.~A.}\ \bibnamefont
  {Bees}},\ }\bibfield  {title} {\enquote {\bibinfo {title} {A tale of three
  taxes: photo-gyro-gravitactic bioconvection},}\ }\href@noop {} {\bibfield
  {journal} {\bibinfo  {journal} {Journal of Experimental Biology}\ }\textbf
  {\bibinfo {volume} {214}},\ \bibinfo {pages} {2398--2408} (\bibinfo {year}
  {2011})}\BibitemShut {NoStop}%
\bibitem [{\citenamefont {Kessler}(1989)}]{5kessler1989}%
  \BibitemOpen
  \bibfield  {author} {\bibinfo {author} {\bibfnamefont {J.}~\bibnamefont
  {Kessler}},\ }\bibfield  {title} {\enquote {\bibinfo {title} {Path and
  pattern-the mutual dynamics of swimming cells and their environment},}\
  }\href@noop {} {\bibfield  {journal} {\bibinfo  {journal} {Comments Theor.
  Biol.}\ }\textbf {\bibinfo {volume} {1}},\ \bibinfo {pages} {85--108}
  (\bibinfo {year} {1989})}\BibitemShut {NoStop}%
\bibitem [{\citenamefont {Ghorai}, \citenamefont {Panda},\ and\ \citenamefont
  {Hill}(2010)}]{7ghorai2010}%
  \BibitemOpen
  \bibfield  {author} {\bibinfo {author} {\bibfnamefont {S.}~\bibnamefont
  {Ghorai}}, \bibinfo {author} {\bibfnamefont {M.~K.}\ \bibnamefont {Panda}}, \
  and\ \bibinfo {author} {\bibfnamefont {N.~A.}\ \bibnamefont {Hill}},\
  }\bibfield  {title} {\enquote {\bibinfo {title} {Bioconvection in a
  suspension of isotropically scattering phototactic algae},}\ }\href@noop {}
  {\bibfield  {journal} {\bibinfo  {journal} {Physics of Fluids}\ }\textbf
  {\bibinfo {volume} {22}},\ \bibinfo {pages} {071901} (\bibinfo {year}
  {2010})}\BibitemShut {NoStop}%
\bibitem [{\citenamefont {Straughan}(1993)}]{9straughan1993}%
  \BibitemOpen
  \bibfield  {author} {\bibinfo {author} {\bibfnamefont {B.}~\bibnamefont
  {Straughan}},\ }\href@noop {} {\emph {\bibinfo {title} {Mathematical aspects
  of penetrative convection}}}\ (\bibinfo  {publisher} {CRC Press},\ \bibinfo
  {year} {1993})\BibitemShut {NoStop}%
\bibitem [{\citenamefont {Ghorai}\ and\ \citenamefont
  {Hill}(2005)}]{10ghorai2005}%
  \BibitemOpen
  \bibfield  {author} {\bibinfo {author} {\bibfnamefont {S.}~\bibnamefont
  {Ghorai}}\ and\ \bibinfo {author} {\bibfnamefont {N.~A.}\ \bibnamefont
  {Hill}},\ }\bibfield  {title} {\enquote {\bibinfo {title} {Penetrative
  phototactic bioconvection},}\ }\href@noop {} {\bibfield  {journal} {\bibinfo
  {journal} {Physics of fluids}\ }\textbf {\bibinfo {volume} {17}},\ \bibinfo
  {pages} {074101} (\bibinfo {year} {2005})}\BibitemShut {NoStop}%
\bibitem [{\citenamefont {Panda}\ and\ \citenamefont
  {Singh}(2016)}]{11panda2016}%
  \BibitemOpen
  \bibfield  {author} {\bibinfo {author} {\bibfnamefont {M.~K.}\ \bibnamefont
  {Panda}}\ and\ \bibinfo {author} {\bibfnamefont {R.}~\bibnamefont {Singh}},\
  }\bibfield  {title} {\enquote {\bibinfo {title} {Penetrative phototactic
  bioconvection in a two-dimensional non-scattering suspension},}\ }\href@noop
  {} {\bibfield  {journal} {\bibinfo  {journal} {Physics of Fluids}\ }\textbf
  {\bibinfo {volume} {28}},\ \bibinfo {pages} {054105} (\bibinfo {year}
  {2016})}\BibitemShut {NoStop}%
\bibitem [{\citenamefont {Vincent}\ and\ \citenamefont
  {Hill}(1996)}]{12vincent1996}%
  \BibitemOpen
  \bibfield  {author} {\bibinfo {author} {\bibfnamefont {R.~V.}\ \bibnamefont
  {Vincent}}\ and\ \bibinfo {author} {\bibfnamefont {N.~A.}\ \bibnamefont
  {Hill}},\ }\bibfield  {title} {\enquote {\bibinfo {title} {Bioconvection in a
  suspension of phototactic algae},}\ }\href@noop {} {\bibfield  {journal}
  {\bibinfo  {journal} {Journal of Fluid Mechanics}\ }\textbf {\bibinfo
  {volume} {327}},\ \bibinfo {pages} {343--371} (\bibinfo {year}
  {1996})}\BibitemShut {NoStop}%
\bibitem [{\citenamefont {Panda}\ and\ \citenamefont
  {Ghorai}(2013)}]{14panda2013}%
  \BibitemOpen
  \bibfield  {author} {\bibinfo {author} {\bibfnamefont {M.~K.}\ \bibnamefont
  {Panda}}\ and\ \bibinfo {author} {\bibfnamefont {S.}~\bibnamefont {Ghorai}},\
  }\bibfield  {title} {\enquote {\bibinfo {title} {Penetrative phototactic
  bioconvection in an isotropic scattering suspension},}\ }\href@noop {}
  {\bibfield  {journal} {\bibinfo  {journal} {Physics of Fluids}\ }\textbf
  {\bibinfo {volume} {25}},\ \bibinfo {pages} {071902} (\bibinfo {year}
  {2013})}\BibitemShut {NoStop}%
\bibitem [{\citenamefont {Panda}\ \emph {et~al.}(2016)\citenamefont {Panda},
  \citenamefont {Singh}, \citenamefont {Mishra},\ and\ \citenamefont
  {Mohanty}}]{15panda2016}%
  \BibitemOpen
  \bibfield  {author} {\bibinfo {author} {\bibfnamefont {M.~K.}\ \bibnamefont
  {Panda}}, \bibinfo {author} {\bibfnamefont {R.}~\bibnamefont {Singh}},
  \bibinfo {author} {\bibfnamefont {A.~C.}\ \bibnamefont {Mishra}}, \ and\
  \bibinfo {author} {\bibfnamefont {S.~K.}\ \bibnamefont {Mohanty}},\
  }\bibfield  {title} {\enquote {\bibinfo {title} {Effects of both diffuse and
  collimated incident radiation on phototactic bioconvection},}\ }\href@noop {}
  {\bibfield  {journal} {\bibinfo  {journal} {Physics of Fluids}\ }\textbf
  {\bibinfo {volume} {28}},\ \bibinfo {pages} {124104} (\bibinfo {year}
  {2016})}\BibitemShut {NoStop}%
\bibitem [{\citenamefont {Panda}(2020)}]{8panda2020}%
  \BibitemOpen
  \bibfield  {author} {\bibinfo {author} {\bibfnamefont {M.~K.}\ \bibnamefont
  {Panda}},\ }\bibfield  {title} {\enquote {\bibinfo {title} {Effects of
  anisotropic scattering on the onset of phototactic bioconvection with diffuse
  and collimated irradiation},}\ }\href@noop {} {\bibfield  {journal} {\bibinfo
   {journal} {Physics of Fluids}\ }\textbf {\bibinfo {volume} {32}},\ \bibinfo
  {pages} {091903} (\bibinfo {year} {2020})}\BibitemShut {NoStop}%
\bibitem [{\citenamefont {Panda}, \citenamefont {Sharma},\ and\ \citenamefont
  {Kumar}(2022)}]{16panda2022}%
  \BibitemOpen
  \bibfield  {author} {\bibinfo {author} {\bibfnamefont {M.~K.}\ \bibnamefont
  {Panda}}, \bibinfo {author} {\bibfnamefont {P.}~\bibnamefont {Sharma}}, \
  and\ \bibinfo {author} {\bibfnamefont {S.}~\bibnamefont {Kumar}},\ }\bibfield
   {title} {\enquote {\bibinfo {title} {Effects of oblique irradiation on the
  onset of phototactic bioconvection},}\ }\href@noop {} {\bibfield  {journal}
  {\bibinfo  {journal} {Physics of Fluids}\ }\textbf {\bibinfo {volume} {34}},\
  \bibinfo {pages} {024108} (\bibinfo {year} {2022})}\BibitemShut {NoStop}%
\bibitem [{\citenamefont {Kumar}(2022)}]{17kumar2022}%
  \BibitemOpen
  \bibfield  {author} {\bibinfo {author} {\bibfnamefont {S.}~\bibnamefont
  {Kumar}},\ }\bibfield  {title} {\enquote {\bibinfo {title} {Phototactic
  isotropic scattering bioconvection with oblique irradiation},}\ }\href@noop
  {} {\bibfield  {journal} {\bibinfo  {journal} {Physics of Fluids}\ }\textbf
  {\bibinfo {volume} {34}},\ \bibinfo {pages} {114125} (\bibinfo {year}
  {2022})}\BibitemShut {NoStop}%
\bibitem [{\citenamefont {Kumar}(2023{\natexlab{a}})}]{39kumar2023}%
  \BibitemOpen
  \bibfield  {author} {\bibinfo {author} {\bibfnamefont {S.}~\bibnamefont
  {Kumar}},\ }\bibfield  {title} {\enquote {\bibinfo {title} {Isotropic
  scattering with a rigid upper surface at the onset of phototactic
  bioconvection},}\ }\href@noop {} {\bibfield  {journal} {\bibinfo  {journal}
  {Physics of Fluids}\ }\textbf {\bibinfo {volume} {35}},\ \bibinfo {pages}
  {024106} (\bibinfo {year} {2023}{\natexlab{a}})}\BibitemShut {NoStop}%
\bibitem [{\citenamefont {Kumar}(2023{\natexlab{b}})}]{40kumar2023}%
  \BibitemOpen
  \bibfield  {author} {\bibinfo {author} {\bibfnamefont {S.}~\bibnamefont
  {Kumar}},\ }\bibfield  {title} {\enquote {\bibinfo {title} {Effect of
  rotation on the suspension of phototactic bioconvection},}\ }\href@noop {}
  {\bibfield  {journal} {\bibinfo  {journal} {Physics of Fluids}\ }\textbf
  {\bibinfo {volume} {35}} (\bibinfo {year} {2023}{\natexlab{b}})}\BibitemShut
  {NoStop}%
\bibitem [{\citenamefont {Panda}\ and\ \citenamefont
  {Rajput}(2023)}]{41rajput2023}%
  \BibitemOpen
  \bibfield  {author} {\bibinfo {author} {\bibfnamefont {M.~K.}\ \bibnamefont
  {Panda}}\ and\ \bibinfo {author} {\bibfnamefont {S.~K.}\ \bibnamefont
  {Rajput}},\ }\bibfield  {title} {\enquote {\bibinfo {title} {Phototactic
  bioconvection with the combined effect of diffuse and oblique collimated flux
  on an algal suspension},}\ }\href@noop {} {\bibfield  {journal} {\bibinfo
  {journal} {Physics of Fluids}\ }\textbf {\bibinfo {volume} {35}} (\bibinfo
  {year} {2023})}\BibitemShut {NoStop}%
\bibitem [{\citenamefont {Cash}\ and\ \citenamefont
  {Moore}(1980)}]{19cash1980}%
  \BibitemOpen
  \bibfield  {author} {\bibinfo {author} {\bibfnamefont {J.~R.}\ \bibnamefont
  {Cash}}\ and\ \bibinfo {author} {\bibfnamefont {D.~R.}\ \bibnamefont
  {Moore}},\ }\bibfield  {title} {\enquote {\bibinfo {title} {A high order
  method for the numerical solution of two-point boundary value problems},}\
  }\href@noop {} {\bibfield  {journal} {\bibinfo  {journal} {BIT Numerical
  Mathematics}\ }\textbf {\bibinfo {volume} {20}},\ \bibinfo {pages} {44--52}
  (\bibinfo {year} {1980})}\BibitemShut {NoStop}%
\bibitem [{\citenamefont {H{\"a}der}(1987)}]{6hader1987}%
  \BibitemOpen
  \bibfield  {author} {\bibinfo {author} {\bibfnamefont {D.-P.}\ \bibnamefont
  {H{\"a}der}},\ }\bibfield  {title} {\enquote {\bibinfo {title} {Polarotaxis,
  gravitaxis and vertical phototaxis in the green flagellate, euglena
  gracilis},}\ }\href@noop {} {\bibfield  {journal} {\bibinfo  {journal}
  {Archives of microbiology}\ }\textbf {\bibinfo {volume} {147}},\ \bibinfo
  {pages} {179--183} (\bibinfo {year} {1987})}\BibitemShut {NoStop}%
\bibitem [{\citenamefont {Ghorai}\ and\ \citenamefont
  {Panda}(2013)}]{13ghorai2013}%
  \BibitemOpen
  \bibfield  {author} {\bibinfo {author} {\bibfnamefont {S.}~\bibnamefont
  {Ghorai}}\ and\ \bibinfo {author} {\bibfnamefont {M.~K.}\ \bibnamefont
  {Panda}},\ }\bibfield  {title} {\enquote {\bibinfo {title} {Bioconvection in
  an anisotropic scattering suspension of phototactic algae},}\ }\href@noop {}
  {\bibfield  {journal} {\bibinfo  {journal} {European Journal of
  Mechanics-B/Fluids}\ }\textbf {\bibinfo {volume} {41}},\ \bibinfo {pages}
  {81--93} (\bibinfo {year} {2013})}\BibitemShut {NoStop}%
\bibitem [{\citenamefont {Hill}\ and\ \citenamefont
  {H{\"a}der}(1997)}]{18hill1997}%
  \BibitemOpen
  \bibfield  {author} {\bibinfo {author} {\bibfnamefont {N.~A.}\ \bibnamefont
  {Hill}}\ and\ \bibinfo {author} {\bibfnamefont {D.-P.}\ \bibnamefont
  {H{\"a}der}},\ }\bibfield  {title} {\enquote {\bibinfo {title} {A biased
  random walk model for the trajectories of swimming micro-organisms},}\
  }\href@noop {} {\bibfield  {journal} {\bibinfo  {journal} {Journal of
  theoretical biology}\ }\textbf {\bibinfo {volume} {186}},\ \bibinfo {pages}
  {503--526} (\bibinfo {year} {1997})}\BibitemShut {NoStop}%
\bibitem [{\citenamefont {Kessler}(1986)}]{25kessler1986}%
  \BibitemOpen
  \bibfield  {author} {\bibinfo {author} {\bibfnamefont {J.~O.}\ \bibnamefont
  {Kessler}},\ }\bibfield  {title} {\enquote {\bibinfo {title} {The external
  dynamics of swimming micro-organisms},}\ }\href@noop {} {\bibfield  {journal}
  {\bibinfo  {journal} {Progress in phycological research}\ }\textbf {\bibinfo
  {volume} {4}},\ \bibinfo {pages} {258--307} (\bibinfo {year}
  {1986})}\BibitemShut {NoStop}%
\bibitem [{\citenamefont {Kessler}\ and\ \citenamefont
  {Hill}(1997)}]{26kessler1997}%
  \BibitemOpen
  \bibfield  {author} {\bibinfo {author} {\bibfnamefont {J.~O.}\ \bibnamefont
  {Kessler}}\ and\ \bibinfo {author} {\bibfnamefont {N.~A.}\ \bibnamefont
  {Hill}},\ }\bibfield  {title} {\enquote {\bibinfo {title} {Complementarity of
  physics, biology and geometry in the dynamics of swimming micro-organisms},}\
  }in\ \href@noop {} {\emph {\bibinfo {booktitle} {Physics of biological
  systems}}}\ (\bibinfo  {publisher} {Springer},\ \bibinfo {year} {1997})\ pp.\
  \bibinfo {pages} {325--340}\BibitemShut {NoStop}%
\bibitem [{\citenamefont {Kage}\ \emph {et~al.}(2013)\citenamefont {Kage},
  \citenamefont {Hosoya}, \citenamefont {Baba},\ and\ \citenamefont
  {Mogami}}]{27kage2013}%
  \BibitemOpen
  \bibfield  {author} {\bibinfo {author} {\bibfnamefont {A.}~\bibnamefont
  {Kage}}, \bibinfo {author} {\bibfnamefont {C.}~\bibnamefont {Hosoya}},
  \bibinfo {author} {\bibfnamefont {S.~A.}\ \bibnamefont {Baba}}, \ and\
  \bibinfo {author} {\bibfnamefont {Y.}~\bibnamefont {Mogami}},\ }\bibfield
  {title} {\enquote {\bibinfo {title} {Drastic reorganization of the
  bioconvection pattern of chlamydomonas: quantitative analysis of the pattern
  transition response},}\ }\href@noop {} {\bibfield  {journal} {\bibinfo
  {journal} {Journal of Experimental Biology}\ }\textbf {\bibinfo {volume}
  {216}},\ \bibinfo {pages} {4557--4566} (\bibinfo {year} {2013})}\BibitemShut
  {NoStop}%
\bibitem [{\citenamefont {Mendelson}\ and\ \citenamefont
  {Lega}(1998)}]{28mendelson1998}%
  \BibitemOpen
  \bibfield  {author} {\bibinfo {author} {\bibfnamefont {N.~H.}\ \bibnamefont
  {Mendelson}}\ and\ \bibinfo {author} {\bibfnamefont {J.}~\bibnamefont
  {Lega}},\ }\bibfield  {title} {\enquote {\bibinfo {title} {A complex pattern
  of traveling stripes is produced by swimming cells of bacillus subtilis},}\
  }\href@noop {} {\bibfield  {journal} {\bibinfo  {journal} {Journal of
  bacteriology}\ }\textbf {\bibinfo {volume} {180}},\ \bibinfo {pages}
  {3285--3294} (\bibinfo {year} {1998})}\BibitemShut {NoStop}%
\bibitem [{\citenamefont {Gittleson}\ and\ \citenamefont
  {Jahn}(1968)}]{29gittleson1968}%
  \BibitemOpen
  \bibfield  {author} {\bibinfo {author} {\bibfnamefont {S.~M.}\ \bibnamefont
  {Gittleson}}\ and\ \bibinfo {author} {\bibfnamefont {T.~L.}\ \bibnamefont
  {Jahn}},\ }\bibfield  {title} {\enquote {\bibinfo {title} {Pattern swimming
  by polytomella agilis},}\ }\href@noop {} {\bibfield  {journal} {\bibinfo
  {journal} {The American Naturalist}\ }\textbf {\bibinfo {volume} {102}},\
  \bibinfo {pages} {413--425} (\bibinfo {year} {1968})}\BibitemShut {NoStop}%
\bibitem [{\citenamefont {Khan}\ \emph {et~al.}(2017)\citenamefont {Khan},
  \citenamefont {Gul}, \citenamefont {Khan}, \citenamefont {Bonyah},\ and\
  \citenamefont {Islam}}]{30khan2017}%
  \BibitemOpen
  \bibfield  {author} {\bibinfo {author} {\bibfnamefont {N.~S.}\ \bibnamefont
  {Khan}}, \bibinfo {author} {\bibfnamefont {T.}~\bibnamefont {Gul}}, \bibinfo
  {author} {\bibfnamefont {M.~A.}\ \bibnamefont {Khan}}, \bibinfo {author}
  {\bibfnamefont {E.}~\bibnamefont {Bonyah}}, \ and\ \bibinfo {author}
  {\bibfnamefont {S.}~\bibnamefont {Islam}},\ }\bibfield  {title} {\enquote
  {\bibinfo {title} {Mixed convection in gravity-driven thin film non-newtonian
  nanofluids flow with gyrotactic microorganisms},}\ }\href@noop {} {\bibfield
  {journal} {\bibinfo  {journal} {Results in physics}\ }\textbf {\bibinfo
  {volume} {7}},\ \bibinfo {pages} {4033--4049} (\bibinfo {year}
  {2017})}\BibitemShut {NoStop}%
\bibitem [{\citenamefont {Hayat}, \citenamefont {Alsaedi}\ \emph
  {et~al.}(2021)\citenamefont {Hayat}, \citenamefont {Alsaedi} \emph
  {et~al.}}]{31hayat2021}%
  \BibitemOpen
  \bibfield  {author} {\bibinfo {author} {\bibfnamefont {T.}~\bibnamefont
  {Hayat}}, \bibinfo {author} {\bibfnamefont {A.}~\bibnamefont {Alsaedi}},
  \emph {et~al.},\ }\bibfield  {title} {\enquote {\bibinfo {title} {Development
  of bioconvection flow of nanomaterial with melting effects},}\ }\href@noop {}
  {\bibfield  {journal} {\bibinfo  {journal} {Chaos, Solitons \& Fractals}\
  }\textbf {\bibinfo {volume} {148}},\ \bibinfo {pages} {111015} (\bibinfo
  {year} {2021})}\BibitemShut {NoStop}%
\bibitem [{\citenamefont {Incropera}, \citenamefont {Wagner},\ and\
  \citenamefont {Houf}(1981)}]{32incropera1981}%
  \BibitemOpen
  \bibfield  {author} {\bibinfo {author} {\bibfnamefont {F.~P.}\ \bibnamefont
  {Incropera}}, \bibinfo {author} {\bibfnamefont {T.~R.}\ \bibnamefont
  {Wagner}}, \ and\ \bibinfo {author} {\bibfnamefont {W.~G.}\ \bibnamefont
  {Houf}},\ }\bibfield  {title} {\enquote {\bibinfo {title} {A comparison of
  predictions and measurements of the radiation field in a shallow water
  layer},}\ }\href@noop {} {\bibfield  {journal} {\bibinfo  {journal} {Water
  Resources Research}\ }\textbf {\bibinfo {volume} {17}},\ \bibinfo {pages}
  {142--148} (\bibinfo {year} {1981})}\BibitemShut {NoStop}%
\bibitem [{\citenamefont {Daniel}, \citenamefont {Laurendeau},\ and\
  \citenamefont {Incropera}(1979)}]{33daniel1979}%
  \BibitemOpen
  \bibfield  {author} {\bibinfo {author} {\bibfnamefont {K.~J.}\ \bibnamefont
  {Daniel}}, \bibinfo {author} {\bibfnamefont {N.~M.}\ \bibnamefont
  {Laurendeau}}, \ and\ \bibinfo {author} {\bibfnamefont {F.~P.}\ \bibnamefont
  {Incropera}},\ }\bibfield  {title} {\enquote {\bibinfo {title} {Prediction of
  radiation absorption and scattering in turbid water bodies},}\ }\href@noop {}
  {\  (\bibinfo {year} {1979})}\BibitemShut {NoStop}%
\bibitem [{\citenamefont {Hill}, \citenamefont {Pedley},\ and\ \citenamefont
  {Kessler}(1989)}]{34hill1989}%
  \BibitemOpen
  \bibfield  {author} {\bibinfo {author} {\bibfnamefont {N.~A.}\ \bibnamefont
  {Hill}}, \bibinfo {author} {\bibfnamefont {T.~J.}\ \bibnamefont {Pedley}}, \
  and\ \bibinfo {author} {\bibfnamefont {J.~O.}\ \bibnamefont {Kessler}},\
  }\bibfield  {title} {\enquote {\bibinfo {title} {Growth of bioconvection
  patterns in a suspension of gyrotactic micro-organisms in a layer of finite
  depth},}\ }\href@noop {} {\bibfield  {journal} {\bibinfo  {journal} {Journal
  of Fluid Mechanics}\ }\textbf {\bibinfo {volume} {208}},\ \bibinfo {pages}
  {509--543} (\bibinfo {year} {1989})}\BibitemShut {NoStop}%
\bibitem [{\citenamefont {Modest}\ and\ \citenamefont
  {Mazumder}(2021)}]{35modest2021}%
  \BibitemOpen
  \bibfield  {author} {\bibinfo {author} {\bibfnamefont {M.~F.}\ \bibnamefont
  {Modest}}\ and\ \bibinfo {author} {\bibfnamefont {S.}~\bibnamefont
  {Mazumder}},\ }\href@noop {} {\emph {\bibinfo {title} {Radiative heat
  transfer}}}\ (\bibinfo  {publisher} {Academic press},\ \bibinfo {year}
  {2021})\BibitemShut {NoStop}%
\bibitem [{\citenamefont {Chandrasekhar}(1960)}]{36chandrasekhar1960}%
  \BibitemOpen
  \bibfield  {author} {\bibinfo {author} {\bibfnamefont {S.}~\bibnamefont
  {Chandrasekhar}},\ }\bibfield  {title} {\enquote {\bibinfo {title} {Radiative
  transfer dover publications inc},}\ }\href@noop {} {\bibfield  {journal}
  {\bibinfo  {journal} {New York}\ } (\bibinfo {year} {1960})}\BibitemShut
  {NoStop}%
\bibitem [{\citenamefont {Ghorai}\ and\ \citenamefont
  {Singh}(2009)}]{37ghorai2009}%
  \BibitemOpen
  \bibfield  {author} {\bibinfo {author} {\bibfnamefont {S.}~\bibnamefont
  {Ghorai}}\ and\ \bibinfo {author} {\bibfnamefont {R.}~\bibnamefont {Singh}},\
  }\bibfield  {title} {\enquote {\bibinfo {title} {Linear stability analysis of
  gyrotactic plumes},}\ }\href@noop {} {\bibfield  {journal} {\bibinfo
  {journal} {Physics of Fluids}\ }\textbf {\bibinfo {volume} {21}},\ \bibinfo
  {pages} {081901} (\bibinfo {year} {2009})}\BibitemShut {NoStop}%
\bibitem [{\citenamefont {Press}(1992)}]{38press1992}%
  \BibitemOpen
  \bibfield  {author} {\bibinfo {author} {\bibfnamefont {W.~H.}\ \bibnamefont
  {Press}},\ }\bibfield  {title} {\enquote {\bibinfo {title} {Numerical recipes
  in fortran.}}\ }\href@noop {} {\bibfield  {journal} {\bibinfo  {journal} {The
  Art of Scientific Computing.}\ } (\bibinfo {year} {1992})}\BibitemShut
  {NoStop}%
\end{thebibliography}%
	
\end{document}